\documentclass[12pt]{article}
\usepackage[]{amsmath,amssymb}
\usepackage{amscd}
\usepackage{latexsym}
\usepackage{cite}

\newtheorem{definition}{Definition}[section]
\newtheorem{theorem}[definition]{Theorem}
\newtheorem{lemma}[definition]{Lemma}
\newtheorem{corollary}[definition]{Corollary}

\newtheorem{note}[definition]{Note}

\newtheorem{proposition}[definition]{Proposition}
\newtheorem{notation}[definition]{Notation}

\def\I{\mathbb I}

\def\C{\mathbb C}

\def\K{\mathbb K}

\begin{document}

\title{Totally Bipartite/ABipartite Leonard pairs and\\Leonard triples of Bannai/Ito type}
\author{George M. F. Brown}
\maketitle

\begin{abstract}
This paper is about three classes of objects: Leonard pairs, Leonard triples, and the finite-dimensional irreducible modules for a certain algebra $\mathcal{A}$.  Let $\K$ denote an algebraically closed field of characteristic zero.  Let $V$ denote a vector space over $\K$ with finite positive dimension.  A Leonard pair on $V$ is an ordered pair of linear transformations in $\mathrm{End}(V)$ such that for each of these transformations there exists a basis for $V$ with respect to which the matrix representing that transformation is diagonal and the matrix representing the other transformation is irreducible tridiagonal.  Whenever these tridiagonal matrices are bipartite, the Leonard pair is said to be totally bipartite.  A mild weakening of the bipartite assumption yields a type of Leonard pair said to be totally almost bipartite.  A Leonard pair is said to be totally B/AB whenever it is totally bipartite or totally almost bipartite.  The notion of a Leonard triple and the corresponding notion of totally B/AB are similarly defined.  There are families of Leonard pairs and Leonard triples said to have Bannai/Ito type.  The Leonard pairs and Leonard triples of interest to us are the ones that are totally B/AB and of Bannai/Ito type.

Let $\mathcal{A}$ denote the unital associative $\K$-algebra defined by generators $x,y,z$ and relations
\[
xy+yx=2z,\qquad yz+zy=2x,\qquad zx+xz=2y.
\]
The algebra $\mathcal{A}$ has a presentation involving generators $x,y$ and relations
\[
x^{2}y+2xyx+yx^{2}=4y,\qquad y^{2}x+2yxy+xy^{2}=4x.
\]

In this paper we obtain the following results.  We classify up to isomorphism the totally B/AB Leonard pairs of Bannai/Ito type.  We classify up to isomorphism the totally B/AB Leonard triples of Bannai/Ito type.  We classify up to isomorphism the finite-dimensional irreducible $\mathcal{A}$-modules.  We show that these three classes of objects are essentially in one-to-one correspondence, and describe these correspondences in detail.
\end{abstract}

\section{Introduction}\label{S:intro}

Throughout this paper, $\K$ denotes an algebraically closed field of characteristic zero.

We now recall the definition of a Leonard pair.  To do this, we use the following terms.  A square matrix $B$ is said to be \emph{tridiagonal} whenever each nonzero entry lies on either the diagonal, the subdiagonal, or the superdiagonal.  Assume $B$ is tridiagonal.  Then $B$ is said to be \emph{irreducible} whenever each entry on the subdiagonal or superdiagonal is nonzero.

\begin{definition}\label{D:LP}
{\normalfont\cite[Definition 1.1]{terwilliger} Let $V$ denote a vector space over $\K$ with finite positive dimension.  By a \emph{Leonard pair} on $V$ we mean an ordered pair of linear transformations $A:V\to V$, $A^{*}:V\to V$ which satisfy the conditions (i), (ii) below.
\begin{itemize}
\item[\rm (i)] There exists a basis for $V$ with respect to which the matrix representing $A$ is diagonal and the matrix representing $A^{*}$ is irreducible tridiagonal.
\item[\rm (ii)] There exists a basis for $V$ with respect to which the matrix representing $A^{*}$ is diagonal and the matrix representing $A$ is irreducible tridiagonal.
\end{itemize}
The \emph{diameter} of the Leonard pair $A,A^{*}$ is defined to be one less than the dimension of $V$.}
\end{definition}

If $A,A^{*}$ is a Leonard pair on $V$ then so is $A^{*},A$.

We will be considering two families of Leonard pairs said to be totally bipartite and totally almost bipartite.  Before defining these families, we first review a few concepts.  Let $V$ denote a vector space over $\K$ with finite positive dimension.  By a \emph{decomposition} of $V$ we mean a sequence of one-dimensional subspaces of $V$ whose direct sum is $V$.  For any basis $\{v_{i}\}_{i=0}^{d}$ for $V$, the sequence $\{\K v_{i}\}_{i=0}^{d}$ is a decomposition of $V$; the decomposition $\{\K v_{i}\}_{i=0}^{d}$ is said to \emph{correspond} to the basis $\{v_{i}\}_{i=0}^{d}$.  Given a decomposition $\{V_{i}\}_{i=0}^{d}$ of $V$, for $0\leq i\leq d$ pick $0\ne v_{i}\in V_{i}$.  Then $\{v_{i}\}_{i=0}^{d}$ is a basis for $V$ which corresponds to $\{V_{i}\}_{i=0}^{d}$.

Let $A,A^{*}$ denote a Leonard pair on $V$.  A basis for $V$ is called \emph{standard} whenever it satisfies Definition \ref{D:LP}(i).  Observe that, given a decomposition $\{V_{i}\}_{i=0}^{d}$ of $V$, the following (i), (ii) are equivalent.
\begin{itemize}
\item[\rm (i)] There exists a standard basis for $V$ which corresponds to $\{V_{i}\}_{i=0}^{d}$.
\item[\rm (ii)] Every basis for $V$ which corresponds to $\{V_{i}\}_{i=0}^{d}$ is standard.
\end{itemize}
We say that the decomposition $\{V_{i}\}_{i=0}^{d}$ is \emph{standard} whenever (i), (ii) hold.  Observe that if the decomposition $\{V_{i}\}_{i=0}^{d}$ is standard, then so is $\{V_{d-i}\}_{i=0}^{d}$ and no other decomposition of $V$ is standard.

For any nonnegative integer $d$ let $\mathrm{Mat}_{d+1}(\K)$ denote the $\K$-algebra consisting of all $d+1$ by $d+1$ matrices that have entries in $\K$.  We index the rows and columns by $0,1,\ldots,d$.

Let $B\in\mathrm{Mat}_{d+1}(\K)$ be tridiagonal.  We say that $B$ is \emph{bipartite} whenever $B_{ii}=0$ for $0\leq i\leq d$.

\begin{definition}\label{D:bipartiteLP}
{\normalfont A Leonard pair $A,A^{*}$ is said to be \emph{bipartite} whenever the matrix representing $A$ from Definition \ref{D:LP}(ii) is bipartite.  The Leonard pair $A,A^{*}$ is said to be \emph{dual bipartite} whenever the Leonard pair $A^{*},A$ is bipartite.  The Leonard pair $A,A^{*}$ is said to be \emph{totally bipartite} whenever it is bipartite and dual bipartite.}
\end{definition}

Let $B\in\mathrm{Mat}_{d+1}(\K)$ be tridiagonal.  We say that $B$ is \emph{almost bipartite} whenever exactly one of $B_{0,0},B_{d,d}$ is nonzero and $B_{ii}=0$ for $1\leq i\leq d-1$.

\begin{definition}\label{D:almostbipartiteLP}
{\normalfont A Leonard pair $A,A^{*}$ is said to be \emph{almost bipartite} whenever the matrix representing $A$ from Definition \ref{D:LP}(ii) is almost bipartite.  The Leonard pair $A,A^{*}$ is said to be \emph{dual almost bipartite} whenever the Leonard pair $A^{*},A$ is almost bipartite.  The Leonard pair $A,A^{*}$ is said to be \emph{totally almost bipartite} whenever it is almost bipartite and dual almost bipartite.}
\end{definition}

The notion of a Leonard triple was introduced by Brian Curtin in \cite{curtin}.  We recall the definition.

\begin{definition}\label{D:LT}
{\normalfont\cite[Definition 1.2]{curtin} Let $V$ denote a vector space over $\K$ with finite positive dimension.  By a \emph{Leonard triple} on $V$ we mean an ordered triple of linear transformations $A:V\to V$, $A^{*}:V\to V$, $A^{\varepsilon}:V\to V$ which satisfy the conditions (i)--(iii) below.
\begin{itemize}
\item[\rm (i)] There exists a basis for $V$ with respect to which the matrix representing $A$ is diagonal and the matrices representing $A^{*}$ and $A^{\varepsilon}$ are irreducible tridiagonal.
\item[\rm (ii)] There exists a basis for $V$ with respect to which the matrix representing $A^{*}$ is diagonal and the matrices representing $A^{\varepsilon}$ and $A$ are irreducible tridiagonal.
\item[\rm (iii)] There exists a basis for $V$ with respect to which the matrix representing $A^{\varepsilon}$ is diagonal and the matrices representing $A$ and $A^{*}$ are irreducible tridiagonal.
\end{itemize}
The \emph{diameter} of the Leonard triple $A,A^{*},A^{\varepsilon}$ is defined to be one less than the dimension of $V$.}
\end{definition}

\begin{definition}\label{D:bipartiteLT}
{\normalfont In Definition \ref{D:LT} we defined a Leonard triple $A,A^{*},A^{\varepsilon}$.  In that definition we mentioned six tridiagonal matrices.  The Leonard triple $A,A^{*},A^{\varepsilon}$ is said to be \emph{totally bipartite} (resp. \emph{totally almost bipartite}) whenever each of the six tridiagonal matrices is bipartite (resp. almost bipartite).}
\end{definition}

For notational convenience, we say that a Leonard pair or Leonard triple is totally B/AB whenever it is either totally bipartite or totally almost bipartite.

For any Leonard triple, any two of the three form a Leonard pair.  We say that these Leonard pairs are \emph{associated} with the Leonard triple.  The Leonard triple is totally bipartite if and only if all of the associated Leonard pairs are totally bipartite.  The Leonard triple is totally almost bipartite if and only if all of the associated Leonard pairs are totally almost bipartite.

In \cite{terwilliger}, Terwilliger classified the Leonard pairs up to isomorphism.  By that classification, the isomorphism classes of Leonard pairs fall naturally into thirteen families: $q$-Racah, $q$-Hahn, dual $q$-Hahn, $q$-Krawtchouk, dual $q$-Krawtchouk, affine $q$-Krawtchouk, quantum $q$-Krawtchouk, Racah, Hahn, dual Hahn, Krawtchouk, Bannai/Ito and orphan.  For each integer $d\geq3$ these families partition the isomorphism classes of Leonard pairs that have diameter $d$.  It remains an open problem to classify the Leonard triples up to isomorphism.  However, in \cite{curtin}, Curtin classified a family of Leonard triples said to be \emph{modular}.

We say that a Leonard triple is of Bannai/Ito type whenever all of its associated Leonard pairs are of Bannai/Ito type.  Leonard pairs of Bannai/Ito type arise in conjunction with the Bannai/Ito polynomials.  These polynomials were introduced in \cite[pp. 271--273]{bannai-ito} by Bannai and Ito.  In \cite{TVZ}, Tsujimoto, Vinet and Zhedanov studied the Bannai/Ito polynomials in conjunction with Dunkl shift operators and representations of Jordan algebras.  Totally B/AB Leonard pairs and Leonard triples also appear in the literature.  In \cite{miklavic}, Miklavi\v{c} studied totally bipartite Leonard triples associated with some representations of the Lie algebra $\mathfrak{sl}_{2}$ constructed using hypercubes.  The Leonard pairs associated with these Leonard triples are of Krawtchouk type.  In \cite{HKP}, Havl\'{i}\v{c}ek, Klimyk and Po\v{s}ta displayed representations of the nonstandard $q$-deformed cyclically symmetric algebra $U'_{q}(\mathfrak{so}_{3})$.  These representations yield both totally bipartite and totally almost bipartite Leonard triples.  The Leonard pairs associated with these Leonard triples are of $q$-Racah type.

The Leonard pairs and Leonard triples of interest to us are the ones that are totally B/AB and of Bannai/Ito type.  To describe these Leonard pairs and Leonard triples, we consider a $\K$-algebra $\mathcal{A}$ defined by generators $x,y,z$ and relations
\begin{equation}\label{E:beginning}
xy+yx=2z,\qquad yz+zy=2x,\qquad zx+xz=2y.
\end{equation}
The algebra $\mathcal{A}$ has an alternate presentation using generators $x,y$ and relations
\[
x^{2}y+2xyx+yx^{2}=4y,\qquad y^{2}x+2yxy+xy^{2}=4x.
\]
The algebra $\mathcal{A}$ has appeared previously in the literature \cite{arik-kayser}.  In \cite[Section 1]{arik-kayser}, Arik and Kayserilioglu introduced an algebra involving the relations (\ref{E:beginning}).  They called this the anticommutator spin algebra and studied it in conjunction with fermionic quantum systems and the angular momentum algebra.  We say more about Arik and Kayserilioglu's results after Theorem \ref{T:class}.

The present paper is about how the following are related: (i) Totally B/AB Leonard pairs of Bannai/Ito type; (ii) Totally B/AB Leonard triples of Bannai/Ito type; (iii) Finite-dimensional irreducible $\mathcal{A}$-modules.  We now summarize our main results.  We classify up to isomorphism the totally B/AB Leonard pairs of Bannai/Ito type.  We classify up to isomorphism the totally B/AB Leonard triples of Bannai/Ito type.  We classify up to isomorphism the finite-dimensional irreducible $\mathcal{A}$-modules.  We show that these three classes of objects are essentially in one-to-one correspondence.  The correspondence is described as follows.  Let $V$ denote a finite-dimensional irreducible $\mathcal{A}$-module.  Then the actions of $x,y$ (resp. $x,y,z$) on $V$ form a totally B/AB Leonard pair (resp. Leonard triple) of Bannai/Ito type.  Conversely, let $A,A^{*}$ (resp. $A,A^{*},A^{\varepsilon}$) denote a totally B/AB Leonard pair (resp. Leonard triple) of Bannai/Ito type with diameter at least $3$ and let $V$ denote the underlying vector space.  Then there exists an irreducible $\mathcal{A}$-module structure on $V$ and nonzero scalars $\xi,\xi^{*}$ (resp. $\xi,\xi^{*},\xi^{\varepsilon}$) such that $A,A^{*}$ (resp. $A,A^{*},A^{\varepsilon}$) act on $V$ as $\xi x,\xi^{*}y$ (resp. $\xi x,\xi^{*}y,\xi^{\varepsilon}z$) respectively.


We now summarize our results in greater detail.  We first describe the algebra $\mathcal{A}$.  As part of this description, we display an action of the symmetric group $S_{4}$ on $\mathcal{A}$ as a group of automorphisms.  We then classify up to isomorphism the finite-dimensional irreducible $\mathcal{A}$-modules.  Let $V$ denote a finite-dimensional irreducible $\mathcal{A}$-module.  We describe how twisting $V$ via an element of $S_{4}$ affects the isomorphism class of $V$.  We obtain the eigenvalues and corresponding primitive idempotents for the actions of $x,y,z$ on $V$.  We use twisting via the $S_{4}$-action to simplify the calculations.  We display six bases for $V$.  With respect to each of these bases the matrix representing one of $x,y,z$ is diagonal and the matrices representing the other two are irreducible tridiagonal.  We display the matrices representing the actions of $x,y,z$ on $V$ with respect to each of the six bases.  From this, we show that $x,y$ act on $V$ as a totally B/AB Leonard pair of Bannai/Ito type and $x,y,z$ act on $V$ as a totally B/AB Leonard triple of Bannai/Ito type.

Next we classify up to isomorphism the totally B/AB Leonard pairs of Bannai/Ito type.  To avoid trivialities, we assume the diameter is at least $3$.  To obtain this classification, we use the Askey-Wilson relations for a Leonard pair $A,A^{*}$ described by Terwilliger and Vidunas \cite{terw-vid}.  For the case in which $A,A^{*}$ is totally B/AB and of Bannai/Ito type, we show that the Askey-Wilson relations take the form
\[
A^{2}A^{*}+2AA^{*}A+A^{*}A^{2}=\varrho A^{*},\qquad A^{*2}A+2A^{*}AA^{*}+AA^{*2}=\varrho^{*}A,
\]
where $\varrho,\varrho^{*}\in\K$ are nonzero.  Using these relations, we show that for every totally B/AB Leonard pair $A,A^{*}$ on $V$ of Bannai/Ito type with diameter at least $3$, there exist nonzero scalars $\xi,\xi^{*}\in\K$ and an $\mathcal{A}$-module structure on $V$ such that $A,A^{*}$ act as $\xi x,\xi^{*}y$ respectively.  From the preceding paragraphs, we obtain a correspondence between finite-dimensional irreducible $\mathcal{A}$-modules and totally B/AB Leonard pairs of Bannai/Ito type.  Using this correspondence we obtain our classification of the totally B/AB Leonard pairs of Bannai/Ito type.

Next we classify up to isomorphism the totally B/AB Leonard triples of Bannai/Ito type.  Again we assume the diameter is at least $3$.  To obtain this classification, we use some results of Nomura and Terwilliger \cite{nom-terw} concerning linear transformations that are tridiagonal with respect to both eigenbases of a Leonard pair $A,A^{*}$.  For the case in which $A,A^{*}$ is associated with a totally B/AB Leonard triple $A,A^{*},A^{\varepsilon}$ of Bannai/Ito type, we use these results to show that
\[
\zeta^{\varepsilon}(AA^{*}+A^{*}A)=A^{\varepsilon},\qquad\zeta(A^{*}A^{\varepsilon}+A^{\varepsilon}A^{*})=A,\qquad\zeta^{*}(A^{\varepsilon}A+AA^{\varepsilon})=A^{*},
\]
where $\zeta,\zeta^{*},\zeta^{\varepsilon}\in\K$ are nonzero.  Using these relations, we show that for every totally B/AB Leonard triple $A,A^{*},A^{\varepsilon}$ on $V$ of Bannai/Ito type with diameter at least $3$, there exist nonzero scalars $\xi,\xi^{*},\xi^{\varepsilon}\in\K$ and an $\mathcal{A}$-module structure on $V$ such that $A,A^{*},A^{\varepsilon}$ act as $\xi x,\xi^{*}y,\xi^{\varepsilon}z$ respectively.  From the preceding paragraphs, we obtain a correspondence between finite-dimensional irreducible $\mathcal{A}$-modules and totally B/AB Leonard triples of Bannai/Ito type.  Using this correspondence we obtain our classification of the totally B/AB Leonard triples of Bannai/Ito type.

The paper is organized as follows.  In Section \ref{S:alg}, we define the algebra $\mathcal{A}$ and display an action of $S_{4}$ on $\mathcal{A}$ as a group of automorphisms.  In Section \ref{S:modules} we classify the finite-dimensional irreducible $\mathcal{A}$-modules.  In Section \ref{S:S4}, we show how twisting a finite-dimensional irreducible $\mathcal{A}$-module via an element of $S_{4}$ affects the isomorphism class of that module.  In Section \ref{S:idempotents}, we work out the primitive idempotents and eigenvalues for the actions of the $\mathcal{A}$-generators $x,y,z$ on a finite-dimensional irreducible $\mathcal{A}$-module.  In Section \ref{S:bases} we display six bases for each finite-dimensional irreducible $\mathcal{A}$-module.  In Section \ref{S:matrices} we display the matrices representing $x,y,z$ with respect to these six bases.  We also show that these actions form a totally B/AB Leonard triple of Bannai/Ito type.  In Section \ref{S:LPclassification}, we classify the totally B/AB Leonard pairs of Bannai/Ito type and show how they correspond to finite-dimensional irreducible $\mathcal{A}$-modules.  In Section \ref{S:LTclassification}, we classify the totally B/AB Leonard triples of Bannai/Ito type and show how they correspond to finite-dimensional irreducible $\mathcal{A}$-modules.

\section{The algebra $\mathcal{A}$ and its automorphisms}\label{S:alg}

We now define the $\K$-algebra $\mathcal{A}$.

\begin{definition}{\normalfont\cite[Section 1]{arik-kayser}}\label{D:A}
{\normalfont Let $\mathcal{A}$ denote the unital associative algebra over $\K$ with generators $x,y,z$ and relations}
\begin{align}
xy+yx&=2z,\label{E:rel1}\\
yz+zy&=2x,\label{E:rel2}\\
zx+xz&=2y.\label{E:rel3}
\end{align}
\end{definition}
Note that $\mathcal{A}$ is generated by any two of $x,y,z$.  This yields the following two-generator presentation of $\mathcal{A}$.

\begin{lemma}\label{L:AltAW}
The algebra $\mathcal{A}$ has a presentation involving generators $x,y$ and relations
\begin{align}
x^{2}y+2xyx+yx^{2}=4y,\label{E:xy1}\\
y^{2}x+2yxy+xy^{2}=4x.\label{E:xy2}
\end{align}
\end{lemma}

\noindent {\it Proof:}
Rewrite relations (\ref{E:rel2}), (\ref{E:rel3}) by eliminating $z$ using line (\ref{E:rel1}).
\hfill $\Box$ \\

\begin{lemma}\label{L:fix}
Any algebra automorphism of $\mathcal{A}$ that fixes at least two of $x,y,z$ is the identity.
\end{lemma}

\noindent {\it Proof:}
Since any two of $x,y,z$ generate $\mathcal{A}$, any automorphism that fixes at least two of $x,y,z$ must fix all of $\mathcal{A}$.
\hfill $\Box$ \\

Each permutation of $x,y,z$ extends to a unique algebra automorphism of $\mathcal{A}$; this can be checked using relations (\ref{E:rel1})--(\ref{E:rel3}).  This gives an action of the symmetric group $S_{3}$ on $\mathcal{A}$ as a group of automorphisms.  There are also algebra automorphisms of $\mathcal{A}$ that change the sign of two of $x,y,z$ while preserving the third; this gives an action of the Klein-four group $K_{4}$ on $\mathcal{A}$ as a group of automorphisms.

In a moment we will show how the $S_{3}$ and $K_{4}$ actions interact, but first it will be useful to establish that these actions are faithful.

\begin{definition}\label{D:I}
{\normalfont Let $\I$ denote the set consisting of the symbols $0,x,y,z$.}
\end{definition}

\begin{lemma}\label{L:f0123}
For $n\in\I$ there exists a unique algebra homomorphism $f_{n}:\mathcal{A}\to\K$ satisfying
\begin{center}
\begin{tabular}{|l|c|c|c|}
\hline
$n$ & $f_{n}(x)$ & $f_{n}(y)$ & $f_{n}(z)$ \\ \hline\hline
$0$ & $1$ & $1$ & $1$ \\ \hline
$x$ & $1$ & $-1$ & $-1$ \\ \hline
$y$ & $-1$ & $1$ & $-1$ \\ \hline
$z$ & $-1$ & $-1$ & $1$ \\ \hline
\end{tabular}
\end{center}
Moreover, $f_{n}$ is surjective.
\end{lemma}

\noindent {\it Proof:}
One verifies that $f_{n}$ exists through routine calculation using Definition \ref{D:A}.  Also $f_{n}$ is unique since $\mathcal{A}$ is generated by $x,y,z$.  Observe $f_{n}$ is nonzero and hence surjective.
\hfill $\Box$ \\

\begin{lemma}\label{L:li}
The elements $x,y,z,1$ are linearly independent in the $\K$-vector space $\mathcal{A}$.
\end{lemma}

\noindent {\it Proof:}
Let $a,b,c,d\in\K$ satisfy $ax+by+cz+d=0$.  For each $n\in\I$, we apply $f_{n}$ to this equation and get
\begin{align*}
a+b+c+d=&0,\\
a-b-c+d=&0,\\
-a+b-c+d=&0,\\
-a-b+c+d=&0.
\end{align*}
The coefficient matrix of the above system of equations is non-singular, so the unique solution is $a=b=c=d=0$.  Therefore $x,y,z,1$ are linearly independent.
\hfill $\Box$ \\

\begin{corollary}\label{C:distinct}
$\pm x,\pm y,\pm z$ are mutually distinct elements of $\mathcal{A}$.
\end{corollary}

\noindent {\it Proof:}
Immediate from Lemma \ref{L:li}.
\hfill $\Box$ \\

Recall the $S_{3}$ and $K_{4}$ actions from below Definition \ref{D:A}.

\begin{corollary}\label{C:faithful}
$S_{3}$ and $K_{4}$ act faithfully on $\mathcal{A}$.
\end{corollary}

\noindent {\it Proof:}
By Corollary \ref{C:distinct}, $S_{3}$ and $K_{4}$ act faithfully on the set $\{\pm x,\pm y,\pm z\}$, so they act faithfully on $\mathcal{A}$.
\hfill $\Box$ \\

We remark that, in Section \ref{S:modules}, we will classify up to isomorphism the finite-dimensional irreducible $\mathcal{A}$-modules.  The solutions to this classification include four infinite classes, corresponding to almost bipartite Leonard triples.  The $\mathcal{A}$-modules in these classes are indexed by a nonnegative integer called the diameter.  The $f_{n}$ from  Lemma \ref{L:f0123} come from the $\mathcal{A}$-modules of diameter $0$ in these classes.

\begin{lemma}\label{L:KK4}
Let $\sigma$ denote an automorphism of $\mathcal{A}$ that fixes each of $x,y,z$ up to sign.  Then $\sigma$ must change the sign of an even number of $x,y,z$.
\end{lemma}

\noindent {\it Proof:}
By Lemma \ref{L:fix}, if $\sigma$ fixes any two of $x,y,z$ it must fix all three, so $\sigma$ cannot change the sign of exactly one of $x,y,z$.  Also, $\sigma$ cannot change the sign of all three of $x,y,z$ because, if it did, we could compose it with a non-identity element of $K_{4}$ to get an automorphism that changes the sign of exactly one of $x,y,z$.  The result follows.
\hfill $\Box$ \\

Let $\mathrm{Aut}(\mathcal{A})$ denote the set consisting of all automorphisms of $\mathcal{A}$ and note that $\mathrm{Aut}(\mathcal{A})$ forms a group under composition.  Let $G$ denote the subgroup of $\mathrm{Aut}(\mathcal{A})$ that fixes the set $\{\pm x,\pm y,\pm z\}$.  Let $S$ denote the subgroup of $\mathrm{Aut}(\mathcal{A})$ that fixes the set $\{x,y,z\}$ and let $K$ denote the subgroup of $\mathrm{Aut}(\mathcal{A})$ that fixes each of $x,y,z$ up to sign.  Observe that $S$ and $K$ are both subgroups of $G$.

Corollary \ref{C:faithful} gives an injection of groups $S_{3}\hookrightarrow\mathrm{Aut}(\mathcal{A})$ and, by construction, the image of this injection is $S$.  Similarly, Corollary \ref{C:faithful} gives an injection of groups $K_{4}\hookrightarrow\mathrm{Aut}(\mathcal{A})$.  By Lemma \ref{L:KK4} and the definition of the $K_{4}$-action, the image of this injection is $K$.  Since $S,K\subseteq G$, this gives group injections, $S_{3}\hookrightarrow G$, $K_{4}\hookrightarrow G$ whose images are $S,K$ respectively.

It will turn out that $G$ is isomorphic to $S_{4}$ and that $G$ is a semi-direct product $K\rtimes S$.

\begin{proposition}\label{P:GS4}
$G=K\rtimes S$.
\end{proposition}

\noindent {\it Proof:}
By \cite[Proposition 11.2]{grillet}, it suffices to show $S\cap K=\{1_{G}\}$, $K\triangleleft G$ and $G=KS$.  By construction, $S\cap K=\{1_{G}\}$.  By definition the elements of $G$ permute $\pm x,\pm y,\pm z$.  We define a binary relation $\sim$ on the set $\{\pm x,\pm y,\pm z\}$ such that $u\sim v$ if and only if $u=\pm v$.  Observe that $\sim$ is an equivalence relation.  Moreover, observe that the elements of $G$ permute the three equivalence classes  of $\sim$, resulting in a group homomorphism $\varphi:G\to S_{3}$.  The kernel of this homomorphism is $K$, so $K\triangleleft G$.  Furthermore, the composition $S\hookrightarrow G\underset{\varphi}{\to}S_{3}$ is an isomorphism $S\to S_{3}$, so $G=KS$.  By these comments $G=K\rtimes S$.
\hfill $\Box$ \\

Our next goal is to show that $G$ is isomorphic to $S_{4}$.

\begin{definition}\label{D:h0123}
{\normalfont For $n\in\I$, define $h_{n}\in\mathcal{A}$ as follows:}
\begin{alignat*}{3}
h_{0}&=x+y+z,\qquad\qquad &h_{x}&=x-y-z,\\
h_{y}&=-x+y-z,\qquad\qquad &h_{z}&=-x-y+z.
\end{alignat*}
\end{definition}

\begin{lemma}\label{L:hgen}
We have
\[
x=\frac{h_{0}+h_{x}}{2},\qquad\qquad y=\frac{h_{0}+h_{y}}{2},\qquad\qquad z=\frac{h_{0}+h_{z}}{2}.
\]
Moreover, the algebra $\mathcal{A}$ is generated by $\{h_{n}\}_{n\in\I}$.
\end{lemma}

\noindent {\it Proof:}
Routine.
\hfill $\Box$ \\

Let $\tilde{G}$ denote the group of all permutations of $\I$ and observe $\tilde{G}$ is isomorphic to $S_{4}$.

\begin{proposition}\label{P:S4action}
There exists a group isomorphism $G\to\tilde{G},\sigma\mapsto\tilde{\sigma}$ such that $\sigma(h_{n})=h_{\tilde{\sigma}(n)}$ for all $n\in\I$.
\end{proposition}

\noindent {\it Proof:}
We first show that $G$ fixes the set $\{h_{n}\}_{n\in\I}$.  Since $G$ is generated by $S$ and $K$ it suffices to show that $S$ and $K$ fix $\{h_{n}\}_{n\in\I}$.  We check that this is the case for $S$ by the construction below Lemma \ref{L:KK4}.  We check that this is the case for $K$ by the construction below Lemma \ref{L:KK4} along with Lemma \ref{L:KK4} itself.  Since $G$ fixes the set $\{h_{n}\}_{n\in\I}$, there is a unique group homomorphism $G\to\tilde{G},\sigma\mapsto\tilde{\sigma}$ such that $\sigma(h_{n})=h_{\tilde{\sigma}(n)}$ for all $n\in\I$.  The action of $G$ on $\{h_{n}\}_{n\in\I}$ is faithful in view of Lemma \ref{L:hgen}.  The homomorphism is an isomorphism since each of $G,\tilde{G}$ have cardinality 24.
\hfill $\Box$ \\

\begin{corollary}\label{C:S4}
The group $G$ is isomorphic to $S_{4}$.
\end{corollary}

\noindent {\it Proof:}
$G$ is isomorphic to $\tilde{G}$ by Proposition \ref{P:S4action} and $\tilde{G}$ is isomorphic to $S_{4}$ by construction.
\hfill $\Box$ \\

We just established a group isomorphism $G\to\tilde{G}$.  We have subgroups $S,K\subseteq G$.  We now consider what this isomorphism does to the elements of $S$ and $K$.  To this end, let $\tilde{S}$ denote the subgroup of $\tilde{G}$ consisting of the elements that fix $0$.  Let $\tilde{K}$ denote the unique normal subgroup of $\tilde{G}$ of order 4.  Note that $\tilde{K}$ consists of
\[(0x)(yz),\qquad\qquad(0y)(zx),\qquad\qquad(0z)(xy),\]
together with the identity.

\begin{lemma}\label{L:S3}
With respect to the group isomorphism $G\to\tilde{G}$ from Proposition \ref{P:S4action}, the image of $S$ is $\tilde{S}$.  Moreover, let $\sigma\in S$.  Recall that $\sigma$ permutes the elements $x,y,z$ of $\mathcal{A}$.  Then $\tilde{\sigma}$ permutes the elements $x,y,z$ of $\I$ in the corresponding way.
\end{lemma}

\noindent {\it Proof:}
First we show how $\sigma$ acts on $h_{0}$.
\begin{align*}
\sigma(h_{0})=&\sigma(x+y+z)\\
=&\sigma(x)+\sigma(y)+\sigma(z)\\
=&x+y+z\\
=&h_{0},
\end{align*}
so $\tilde{\sigma}$ fixes $0$.  Let $a,b,c$ denote distinct elements of $\{x,y,z\}$.  Then
\begin{align*}
\sigma(h_{a})=&\sigma(a-b-c)\\
=&\sigma(a)-\sigma(b)-\sigma(c),
\end{align*}
so $\tilde{\sigma}(a)=\sigma(a)$ when $a\in\{x,y,z\}$.  The image of $S$ is $\tilde{S}$ by the definition of $\tilde{S}$.
\hfill $\Box$ \\

\begin{lemma}\label{L:K4}
With respect to the isomorphism $G\to\tilde{G}$ from Proposition \ref{P:S4action}, the image of $K$ is $\tilde{K}$.  Given a non-identity element $\sigma\in K$, recall that $\sigma$ fixes one of $x,y,z$ and changes the sign of the other two.  Let $a,b,c$ denote distinct elements of $\{x,y,z\}$ such that $\sigma$ fixes $a$ and changes the sign of $b$ and $c$.  Now, viewing $a,b,c$ as elements of $\I$, then $\tilde{\sigma}$ is $(0,a)(b,c)$.
\end{lemma}

\noindent {\it Proof:}
$\sigma,\tilde{\sigma}$ are both involutions, so $\tilde{\sigma}$ is a composition of disjoint 2-cycles.  It is therefore sufficient to show how $\tilde{\sigma}$ acts on $0$ and $b$.
\begin{align*}
\sigma(h_{0})=&\sigma(a+b+c)\\
=&\sigma(a)+\sigma(b)+\sigma(c)\\
=&a-b-c\\
=&h_{a},
\end{align*}
so $\tilde{\sigma}$ switches $0$ and $a$.
\begin{align*}
\sigma(h_{b})=&\sigma(-a+b-c)\\
=&-\sigma(a)+\sigma(b)-\sigma(c)\\
=&-a-b+c\\
=&h_{c},
\end{align*}
so $\tilde{\sigma}$ switches $b$ and $c$.  The image of $K$ is $\tilde{K}$ by the definition of $\tilde{K}$.
\hfill $\Box$ \\

\section{The finite-dimensional irreducible $\mathcal{A}$-modules}\label{S:modules}

In this section we classify the finite-dimensional irreducible $\mathcal{A}$-modules up to isomorphism.  This classification is given in Theorem \ref{T:class}.

We adopt the following conventions.  Let $V$ denote a vector space over $\K$.  By $\rm{End}$$(V)$ we mean the $\K$-algebra of linear transformations from $V$ to $V$.  Let $B\in\rm{End}$$(V)$.  By an \emph{eigenvalue} of $B$ we mean a root of the minimal polynomial of $B$.  For an eigenvalue $\theta$ of $B$, the \emph{eigenspace for $B$ associated with $\theta$} is the subspace $\{v\in V|B.v=\theta v\}$.  $B$ is \emph{diagonalizable} whenever $V$ is spanned by its eigenspaces.

\begin{definition}\label{D:vlambda}
{\normalfont Let $V$ denote an $\mathcal{A}$-module.  For $\lambda\in\K$, we define $V(\lambda)=\{v\in V|x.v=\lambda v\}$.}
\end{definition}

\begin{lemma}\label{L:tridi}
Let $V$ denote an $\mathcal{A}$-module.  Then $(y+z).V(\lambda)\subseteq V(2-\lambda)$ and $(y-z).V(\lambda)\subseteq V(-2-\lambda)$.  Moreover, $y.V(\lambda)\subseteq V(2-\lambda)+V(-2-\lambda)$ and $z.V(\lambda)\subseteq V(2-\lambda)+V(-2-\lambda)$.
\end{lemma}

\noindent {\it Proof:}
Let $v\in V(\lambda)$.  Using Definition \ref{D:A} we find that $(y+z).v\in V(2-\lambda)$ and $(y-z).v\in V(-2-\lambda)$.  The first two assertions follow from this.  The last two assertions follow from the first two and the observation that each of $y$ and $z$ is a linear combination of $y+z,y-z$.
\hfill $\Box$ \\

We define functions $f:\K\to\K$ and $g:\K\to\K$ such that $f(\lambda)=2-\lambda$ and $g(\lambda)=-2-\lambda$ for all $\lambda\in\K$.  Observe $f(f(\lambda))=\lambda$ and $g(g(\lambda))=\lambda$ for all $\lambda\in\K$, so $f$ and $g$ are permutations of $\K$.  Note that $f$ has a single orbit of size $1$, namely $\{1\}$ and all other orbits have size $2$.  Similarly, $g$ has a single orbit of size $1$, namely $\{-1\}$ and all other orbits have size $2$.

We make an observation.

\begin{lemma}\label{L:orbits}
The sum of the elements in an orbit of $f$ is equal to the size of the orbit.  The sum of the elements of an orbit of $g$ is equal to $-1$ times the size of the orbit.
\end{lemma}

\begin{definition}\label{D:closed}
{\normalfont Given a set $L$ of elements of $\K$, we say that $L$ is \emph{closed} whenever $f(L)\subseteq L$ and $g(L)\subseteq L$.}
\end{definition}

\begin{lemma}\label{L:clocon}
Let $L$ denote a nonempty closed subset of $\K$.  Then $L$ has infinitely many elements.
\end{lemma}

\noindent {\it Proof:}
We assume $L$ has finite cardinality $n$ and obtain a contradiction.  Because $L$ is closed, it can be partitioned into orbits of $f$.  By Lemma \ref{L:orbits}, the sum of the elements in $L$ is $n$.  Similarly, $L$ can be partitioned into orbits of $g$.  By Lemma \ref{L:orbits}, the sum of the elements in $L$ is $-n$.  This implies $n=-n$, so $n=0$.  But $L$ is nonempty, a contradiction.  The result follows.
\hfill $\Box$ \\

\begin{definition}\label{D:connected}
{\normalfont We say that two distinct elements of $\K$ are \emph{adjacent} whenever they are in the same $f$-orbit or the same $g$-orbit.  A set $L\subseteq\K$ is said to be \emph{connected} whenever the following (i), (ii) hold.
\begin{enumerate}
\item[\rm (i)] $L$ is nonempty.

\item[\rm (ii)] For any partition of $L$ into nonempty subsets $M_{1}$ and $M_{2}$ there exist $\mu\in M_{1}$ and $\sigma\in M_{2}$ such that $\mu$ and $\sigma$ are adjacent.
\end{enumerate}
}
\end{definition}

\begin{lemma}\label{L:span}
Let $V$ denote a finite-dimensional irreducible $\mathcal{A}$-module.  Then the action of $x$ on $V$ is diagonalizable.  Moreover, the set $L=\{\lambda\in\K|V(\lambda)\ne0\}$ is connected.
\end{lemma}

\noindent {\it Proof:}
Since $V$ is nonzero and finite-dimensional and since the ground field $\K$ is algebraically closed there exists a nonzero vector in $V$ that is an eigenvector for $x$.  Therefore $V(\lambda)\ne0$ where $\lambda$ is the corresponding eigenvalue.  So $L$ is nonempty.

Let $M_{1},M_{2}$ denote a partition of $L$ such that $M_{1}$ is nonempty and no element of $M_{1}$ is adjacent to any element of $M_{2}$.  Define $W=\sum_{\mu\in M_{1}}V(\mu)$.  Then $W$ is closed under the action of $\mathcal{A}$ by Lemma \ref{L:tridi}, and nonzero because $M_{1}$ is nonempty and $V(\lambda)\ne 0$ for all $\lambda\in M_{1}$.

Since the $\mathcal{A}$-module $V$ is irreducible, we have $V=W$.  It follows that $M_{1}=L$ and $M_{2}$ is empty, so $L$ is connected.  Furthermore, we have $V=\sum_{\mu\in L}V(\mu)$, so the action of $x$ on $V$ is diagonalizable.
\hfill $\Box$ \\

We will continue discussing the finite-dimensional irreducible $\mathcal{A}$-modules after a comment.

\begin{lemma}\label{L:path}
Let $L$ denote a finite and connected subset of $\K$ with cardinality $d+1$.  Then there is an ordering $\{\theta_{i}\}_{i=0}^{d}$ of the elements of $L$ such that $\theta_{i},\theta_{i+1}$ are adjacent for $0\leq i\leq d-1$.
\end{lemma}

\noindent {\it Proof:}
We will construct an ordering $\{\theta_{i}\}_{i=0}^{d}$ of the elements of $L$.  Assume $d\geq1$; otherwise, the result is trivial.  By definition, $L$ is finite and nonempty.  Therefore, by Lemma \ref{L:clocon}, $L$ is not closed, so there must be an element $\theta_{0}\in L$ such that either $f(\theta_{0})\notin L$ or $g(\theta_{0})\notin L$.  Exactly one of $f(\theta_{0}),g(\theta_{0})$ is in $L$ or else the sets $\{\theta_{0}\}$ and $L\setminus\{\theta_{0}\}$ will violate Definition \ref{D:connected}(ii).  If $f(\theta_{0})\in L$ define $\{\theta_{i}\}_{i=0}^{d}$ to be the first $d+1$ elements of the sequence
\[
\theta_{0},\quad f(\theta_{0}),\quad g(f(\theta_{0})),\quad f(g(f(\theta_{0}))),\quad g(f(g(f(\theta_{0})))),\ldots
\]
If $g(\theta_{0})\in L$ define $\{\theta_{i}\}_{i=0}^{d}$ to be the first $d+1$ elements of the sequence
\[
\theta_{0},\quad g(\theta_{0}),\quad f(g(\theta_{0})),\quad g(f(g(\theta_{0}))),\quad f(g(f(g(\theta_{0})))),\ldots
\]
We claim that $\{\theta_{i}\}_{i=0}^{d}$ is an ordering of the elements of $L$.  Of the integers $0,1,\ldots,d$, let $c$ denote the maximal one such that $\{\theta_{i}\}_{i=0}^{c}$ are mutually distinct and in $L$.  We show that $c=d$.  Let $M_{1}=\{\theta_{i}\}_{i=0}^{c}$ and $M_{2}=L\setminus M_{1}$.  Then $M_{1},M_{2}$ is a partition of $L$ and no element of $M_{1}$ is adjacent to an element of $M_{2}$.  By Definition \ref{D:connected}(ii), one of $M_{1},M_{2}$ is empty.  By construction $M_{1}$ is nonempty so $M_{2}$ is empty and $M_{1}=L$.  Therefore $c=d$, thus proving the claim.  By construction $\theta_{i},\theta_{i+1}$ are adjacent for $0\leq i\leq d-1$.  The result follows.
\hfill $\Box$ \\

\begin{corollary}\label{C:path}
Let $V$ denote a finite-dimensional irreducible $\mathcal{A}$-module.  Then there is an ordering $\{\theta_{i}\}_{i=0}^{d}$ of the eigenvalues for the action of $x$ on $V$ such that $\theta_{i},\theta_{i+1}$ are adjacent for $0\leq i\leq d-1$.
\end{corollary}

\noindent {\it Proof:}
Immediate from Lemmas \ref{L:span} and \ref{L:path}.
\hfill $\Box$ \\

Let $V$ denote a finite-dimensional irreducible $\mathcal{A}$-module.  An ordering $\{\theta_{i}\}_{i=0}^{d}$ of elements of $\K$ will be called \emph{standard} whenever $\theta_{i},\theta_{i+1}$ are adjacent for $0\leq i\leq d-1$.  Note that if the ordering $\{\theta_{i}\}_{i=0}^{d}$ is standard then so is the ordering $\{\theta_{d-i}\}_{i=0}^{d}$.  When we display our Leonard pairs and Leonard triples it will turn out that the eigenvalues for the action of $x$ on a standard decomposition of $V$ form a standard ordering of the eigenvalues.

Let $\{\theta_{i}\}_{i=0}^{d}$ denote a standard ordering of eigenvalues for the action of $x$ on $V$.  For $d\geq 1$,
\begin{equation}\label{E:eigvals}
\theta_{i}=(-1)^{i}(\theta_{0}-2\varepsilon i)\qquad(0\leq i\leq d),
\end{equation}
where $\varepsilon=1$ if $\theta_{1}=f(\theta_{0})$ and $\varepsilon=-1$ if $\theta_{1}=g(\theta_{0})$.  Note that, for $d=0$, equation (\ref{E:eigvals}) holds for $\varepsilon=\pm1$.

We now consider how an element in $\{\theta_{i}\}_{i=0}^{d}$ could be adjacent to an element of $\K$ not among $\{\theta_{i}\}_{i=0}^{d}$.  Recall that if $\lambda,\mu\in\K$ are adjacent then either $\lambda=f(\mu)$ or $\lambda=g(\mu)$.  First assume that $d=0$.  Then $\theta_{0}$ is adjacent to a number other than $\theta_{0}$ because $f(\theta_{0})\ne g(\theta_{0})$.  Next assume that $d\geq1$.  By construction, $\theta_{j}$ is adjacent only to $\theta_{j-1},\theta_{j+1}$ for $1\leq j\leq d-1$.

\begin{lemma}\label{L:one}
With the above notation, assume $d\geq1$.  The following table holds.
\begin{center}
\begin{tabular}{|l|l|l|l|}
\hline
$\varepsilon$ & $\theta_{0}$ & $\mbox{\em{Values for $f$ and $g$}}$ & $\mbox{\em{$\theta_{0}$ is adjacent to}}$ \\ \hline\hline
$1$ & $-1$ & $f(\theta_{0})=\theta_{1},g(\theta_{0})=\theta_{0}$ & only $\theta_{1}$ \\ \hline
$1$ & $\ne-1$ & $f(\theta_{0})=\theta_{1},g(\theta_{0})\notin\{\theta_{i}\}_{i=0}^{d}$ & $\theta_{1}$ and an element of $\K\setminus\{\theta_{i}\}_{i=0}^{d}$ \\ \hline
$-1$ & $1$ & $f(\theta_{0})=\theta_{0},g(\theta_{0})=\theta_{1}$ & only $\theta_{1}$ \\ \hline
$-1$ & $\ne1$ & $f(\theta_{0})\notin\{\theta_{i}\}_{i=0}^{d},g(\theta_{0})=\theta_{1}$ & $\theta_{1}$ and an element of $\K\setminus\{\theta_{i}\}_{i=0}^{d}$ \\ \hline
\end{tabular}
\end{center}
Define $\varepsilon'=(-1)^{d-1}\varepsilon$ and note that $\varepsilon'=1$ if $\theta_{d-1}=f(\theta_{d})$ and $\varepsilon'=-1$ if $\theta_{d-1}=g(\theta_{d})$.  Then the following table holds.
\begin{center}
\begin{tabular}{|l|l|l|l|}
\hline
$\varepsilon'$ & $\theta_{d}$ & $\mbox{\em{Values for $f$ and $g$}}$ & $\mbox{\em{$\theta_{d}$ is adjacent to}}$ \\ \hline\hline
$1$ & $-1$ & $f(\theta_{d})=\theta_{d-1},g(\theta_{d})=\theta_{d}$ & only $\theta_{d-1}$ \\ \hline
$1$ & $\ne-1$ & $f(\theta_{d})=\theta_{d-1},g(\theta_{d})\notin\{\theta_{i}\}_{i=0}^{d}$ & $\theta_{d-1}$ and an element of $\K\setminus\{\theta_{i}\}_{i=0}^{d}$ \\ \hline
$-1$ & $1$ & $f(\theta_{d})=\theta_{d},g(\theta_{d})=\theta_{d-1}$ & only $\theta_{d-1}$ \\ \hline
$-1$ & $\ne1$ & $f(\theta_{d})\notin\{\theta_{i}\}_{i=0}^{d},g(\theta_{d})=\theta_{d-1}$ & $\theta_{d-1}$ and an element of $\K\setminus\{\theta_{i}\}_{i=0}^{d}$ \\ \hline
\end{tabular}
\end{center}
\end{lemma}

\noindent {\it Proof:}
We first show the first table holds.  Rows $1,3$: immediate.

Row 2: By construction $f(\theta_{0})=\theta_{1}$.  We now show that $g(\theta_{0})\notin\{\theta_{i}\}_{i=0}^{d}$.  By way of contradiction, assume $g(\theta_{0})\in\{\theta_{i}\}_{i=0}^{d}$.  Then there exists an integer $i$ with $0\leq i\leq d$ such that $g(\theta_{0})=\theta_{i}$.  By (\ref{E:eigvals}), the definition of $g$, and the fact that $\varepsilon=1$, we have
\begin{equation}\label{E:one}
-2-\theta_{0}=(-1)^{i}(\theta_{0}-2i).
\end{equation}
First assume $i$ is odd.  Then (\ref{E:one}) reduces to $i=-1$, a contradiction.  Next assume $i$ is even.  Then (\ref{E:one}) reduces to $\theta_{0}=i-1$.  We now show that $i=0$.  Assume not.  Then, by (\ref{E:eigvals}) with $i-1$ we find that $\theta_{i-1}=\theta_{0}$ but $i-1\ne0$ since $i$ is even.  This contradicts the fact that $\{\theta_{i}\}_{i=0}^{d}$ are distinct.  Therefore $i=0$ so $\theta_{0}=-1$, a contradiction.  We have now shown that $g(\theta_{0})\notin\{\theta_{i}\}_{i=0}^{d}$.  It follows that $\theta_{0}$ is adjacent to $\theta_{1}$ and an element of $\K\setminus\{\theta_{i}\}_{i=0}^{d}$.

Row $4$: similar to row $2$.

To obtain Table 2, apply Table 1 to the standard ordering $\{\theta_{d-i}\}_{i=0}^{d}$ of eigenvalues for the action of $x$ on $V$.
\hfill $\Box$ \\

We will be discussing five classes of $\mathcal{A}$-modules.  The first class will be denoted $B(d)$ ($B$ for ``bipartite").  The other four will be denoted $AB(d,n)$ with $n\in\I$ ($AB$ for ``almost bipartite").  It will become clear in Section \ref{S:matrices} why we use these terms.  We now introduce the first of these classes.

\begin{lemma}\label{L:Bd}
Let $d$ denote a nonnegative even integer.  There exists an $\mathcal{A}$-module $V$ with basis $\{v_{i}\}_{i=0}^{d}$ on which $x,y,z$ act as follows.  For $0\leq i\leq d$,
\begin{align}
x.v_{i}=&(-1)^{i}(d-2i)v_{i},\label{E:Bd1}\\
y.v_{i}=&(d-i+1)v_{i-1}+(i+1)v_{i+1},\label{E:Bd2}\\
z.v_{i}=&(-1)^{i-1}(d-i+1)v_{i-1}+(-1)^{i}(i+1)v_{i+1},\label{E:Bd3}
\end{align}
where $v_{-1}=0$ and $v_{d+1}=0$.  Moreover $V$ is irreducible.  An $\mathcal{A}$-module of this isomorphism class is said to have type $B(d)$.
\end{lemma}

\noindent {\it Proof:}
One can show that $V$ is an $\mathcal{A}$-module by routine calculation using Definition \ref{D:A}.  We now show that $V$ is irreducible.  Let $W$ denote a nonzero $\mathcal{A}$-submodule of $V$.  We claim that for $0\leq i\leq d-1$, if $v_{i}\in W$ then $v_{i+1}\in W$.  Let $i$ be given and assume $v_{i}\in W$.  Adding (\ref{E:Bd2}) to $(-1)^{i}$ times (\ref{E:Bd3}), we find $(y+(-1)^{i}z).v_{i}=2(i+1)v_{i+1}$.  Because $2(i+1)$ is nonzero, we have $v_{i+1}\in W$ as desired.  A similar argument shows that, for $1\leq i\leq d$, if $v_{i}\in W$ then $v_{i-1}\in W$.

We now show that there exists an integer $j$ $(0\leq j\leq d)$ such that $v_{j}\in W$.  For notational convenience define $\theta_{i}=(-1)^{i}(d-2i)$ for $0\leq i\leq d$ and consider the following elements of $\mathcal{A}$:

\begin{equation}\label{E:exi}
e_{i}=\prod_{
\genfrac{}{}{0pt}{} {0\leq j\leq d}{j\ne i}
}\frac{x-\theta_{j}1}{\theta_{i}-\theta_{j}}\qquad(0\leq i\leq d).
\end{equation}

Using (\ref{E:Bd1}), we obtain $e_{i}.v_{j}=\delta_{ij}v_{j}$ for $0\leq i,j\leq d$.  Here $\delta_{ij}$ denotes the Kronecker delta.  Recall that $\{v_{i}\}_{i=0}^{d}$ is a basis for $V$.  Let $v=c_{0}v_{0}+c_{1}v_{1}+\cdots+c_{d}v_{d}$ denote a nonzero vector in $W$.  Since $v$ is nonzero there exists $j$ $(0\leq j\leq d)$ such that $c_{j}$ is nonzero.  Then $e_{j}.v=c_{j}v_{j}$ is a nonzero scalar multiple of $v_{j}$, so $v_{j}\in W$.  By this and our preliminary comments we find that $W=V$.
\hfill $\Box$ \\

\begin{note}\label{N:irred}
{\normalfont For $d$ odd, an $\mathcal{A}$-module $V$ as in Lemma \ref{L:Bd} exists, but it is not irreducible.  Indeed, we have a direct sum of $\mathcal{A}$-modules $V=V_{1}+V_{2}$ where $V_{1}=\mathrm{span}\{v_{i}+v_{d-i}\}_{i=0}^{d}$ and $V_{2}=\mathrm{span}\{v_{i}-v_{d-i}\}_{i=0}^{d}$.}
\end{note}

\begin{lemma}\label{L:ABd0}
Let $d$ denote a nonnegative integer.  There exists an $\mathcal{A}$-module $V$ with basis $\{v_{i}\}_{i=0}^{d}$ on which $x,y,z$ act as follows.  For $0\leq i\leq d$,
\begin{align}
x.v_{i}=&(-1)^{d+i}(2d-2i+1)v_{i},\label{E:ABd01}\\
y.v_{i}=&(-1)^{d}(2d-i+2)v_{i-1}+(-1)^{d}(i+1)v_{i+1},\label{E:ABd02}\\
z.v_{i}=&(-1)^{i-1}(2d-i+2)v_{i-1}+(-1)^{i}(i+1)v_{i+1},\label{E:ABd03}
\end{align}
where $v_{-1}=0$ and $v_{d+1}=v_{d}$.  Moreover $V$ is irreducible.  An $\mathcal{A}$-module of this isomorphism class is said to have type $AB(d,0)$.
\end{lemma}

\noindent {\it Proof:}
Similar to the proof of Lemma \ref{L:Bd}.
\hfill $\Box$ \\

\begin{lemma}\label{L:ABdx}
Let $d$ denote a nonnegative integer.  There exists an $\mathcal{A}$-module $V$ with basis $\{v_{i}\}_{i=0}^{d}$ on which $x,y,z$ act as follows.  For $0\leq i\leq d$,
\begin{align}
x.v_{i}=&(-1)^{d+i}(2d-2i+1)v_{i},\label{E:ABdx1}\\
y.v_{i}=&(-1)^{d+1}(2d-i+2)v_{i-1}+(-1)^{d+1}(i+1)v_{i+1},\label{E:ABdx2}\\
z.v_{i}=&(-1)^{i}(2d-i+2)v_{i-1}+(-1)^{i+1}(i+1)v_{i+1},\label{E:ABdx3}
\end{align}
where $v_{-1}=0$ and $v_{d+1}=v_{d}$.  Moreover $V$ is irreducible.  An $\mathcal{A}$-module of this isomorphism class is said to have type $AB(d,x)$.
\end{lemma}

\noindent {\it Proof:}
Similar to the proof of Lemma \ref{L:Bd}.
\hfill $\Box$ \\

\begin{lemma}\label{L:ABdy}
Let $d$ denote a nonnegative integer.  There exists an $\mathcal{A}$-module $V$ with basis $\{v_{i}\}_{i=0}^{d}$ on which $x,y,z$ act as follows.  For $0\leq i\leq d$,
\begin{align}
x.v_{i}=&(-1)^{d+i+1}(2d-2i+1)v_{i},\label{E:ABdy1}\\
y.v_{i}=&(-1)^{d}(2d-i+2)v_{i-1}+(-1)^{d}(i+1)v_{i+1},\label{E:ABdy2}\\
z.v_{i}=&(-1)^{i}(2d-i+2)v_{i-1}+(-1)^{i+1}(i+1)v_{i+1},\label{E:ABdy3}
\end{align}
where $v_{-1}=0$ and $v_{d+1}=v_{d}$.  Moreover $V$ is irreducible.  An $\mathcal{A}$-module of this isomorphism class is said to have type $AB(d,y)$.
\end{lemma}

\noindent {\it Proof:}
Similar to the proof of Lemma \ref{L:Bd}.
\hfill $\Box$ \\

\begin{lemma}\label{L:ABdz}
Let $d$ denote a nonnegative integer.  There exists an $\mathcal{A}$-module $V$ with basis $\{v_{i}\}_{i=0}^{d}$ on which $x,y,z$ act as follows.  For $0\leq i\leq d$,
\begin{align}
x.v_{i}=&(-1)^{d+i+1}(2d-2i+1)v_{i},\label{E:ABdz1}\\
y.v_{i}=&(-1)^{d+1}(2d-i+2)v_{i-1}+(-1)^{d+1}(i+1)v_{i+1},\label{E:ABdz2}\\
z.v_{i}=&(-1)^{i-1}(2d-i+2)v_{i-1}+(-1)^{i}(i+1)v_{i+1},\label{E:ABdz3}
\end{align}
where $v_{-1}=0$ and $v_{d+1}=v_{d}$.  Moreover $V$ is irreducible.  An $\mathcal{A}$-module of this isomorphism class is said to have type $AB(d,z)$.
\end{lemma}

\noindent {\it Proof:}
Similar to the proof of Lemma \ref{L:Bd}.
\hfill $\Box$ \\

\begin{definition}\label{D:diameter}
{\normalfont Let $V$ denote a finite-dimensional irreducible $\mathcal{A}$-module from Lemmas \ref{L:Bd}--\ref{L:ABdz}.  We define the \emph{diameter} of $V$ to be one less than the dimension of $V$.  Thus $\mathcal{A}$-modules of types $B(d)$ and $AB(d,n)$ have diameter $d$.}
\end{definition}

\begin{definition}\label{D:BAB}
{\normalfont An $\mathcal{A}$-module $V$ is said to have \emph{type} $B$ when there exists an even integer $d\geq 0$ such that $V$ is of type $B(d)$.  The module is said to have \emph{type} $AB$ when there exists an integer $d\geq 0$ and $n\in\I$ such that $V$ is of type $AB(d,n)$.}
\end{definition}

We comment on Definition \ref{D:BAB}.  We will explain in Section \ref{S:matrices} that on an $\mathcal{A}$-module of type $B$, the generators $x,y,z$ act as a totally bipartite Leonard triple and on an $\mathcal{A}$-module of type $AB$, the generators $x,y,z$ act as a totally almost bipartite Leonard triple.

Our goal for the rest of this section is to show that every finite-dimensional irreducible $\mathcal{A}$-module is isomorphic to exactly one $\mathcal{A}$-module from Lemmas \ref{L:Bd}--\ref{L:ABdz}.  As the next result shows, we can distinguish between the five families using the traces of the $x,y,z$ actions.

\begin{theorem}\label{T:trace}
Let $V$ denote an $\mathcal{A}$-module contained in one of the five families from Lemmas \ref{L:Bd}--\ref{L:ABdz}.  Then the traces of $x,y,z$ on $V$ are given in the following table.

\begin{center}
\begin{tabular}{|l||c|c|c|}
\hline
& $\mbox{\em{tr}}(x)$ & $\mbox{\em{tr}}(y)$ & $\mbox{\em{tr}}(z)$ \\ \hline\hline
$B(d)$ & $0$ & $0$ & $0$ \\ \hline
$AB(d,0)$ & $(-1)^{d}(d+1)$ & $(-1)^{d}(d+1)$ & $(-1)^{d}(d+1)$ \\ \hline
$AB(d,x)$ & $(-1)^{d}(d+1)$ & $(-1)^{d+1}(d+1)$ & $(-1)^{d+1}(d+1)$ \\ \hline
$AB(d,y)$ & $(-1)^{d+1}(d+1)$ & $(-1)^{d}(d+1)$ & $(-1)^{d+1}(d+1)$ \\ \hline
$AB(d,z)$ & $(-1)^{d+1}(d+1)$ & $(-1)^{d+1}(d+1)$ & $(-1)^{d}(d+1)$ \\ \hline
\end{tabular}
\end{center}
\end{theorem}

\noindent {\it Proof:}
Routine.
\hfill $\Box$ \\

\begin{theorem}\label{T:class}
Every finite-dimensional irreducible $\mathcal{A}$-module is isomorphic to exactly one of the modules from Lemmas \ref{L:Bd}--\ref{L:ABdz}.
\end{theorem}

\noindent {\it Proof:}
We first claim that the modules from Lemmas \ref{L:Bd}--\ref{L:ABdz} are mutually non-isomorphic.  To do this we refer to the table from Theorem \ref{T:trace}.  If two such $\mathcal{A}$-modules have different values of $d$, then they have different dimensions and are therefore non-isomorphic.  If they have the same value of $d$, but come from different rows of the table, then they must differ on the traces of at least one of $x,y,z$ and are therefore non-isomorphic.  The claim follows.

Let $V$ denote a finite-dimensional irreducible $\mathcal{A}$-module.  We will show that $V$ is isomorphic to a module from Lemmas \ref{L:Bd}--\ref{L:ABdz}.  Let $\{\theta_{i}\}_{i=0}^{d}$ denote a standard ordering of the eigenvalues for the action of $x$ on $V$.  Recall that the ordering $\{\theta_{d-i}\}_{i=0}^{d}$ is also standard.

By Lemma \ref{L:clocon}, there exists an integer $r$ $(0\leq r\leq d)$ such that $\theta_{r}$ is adjacent to an element of $\K$ not among $\{\theta_{i}\}_{i=0}^{d}$.  By the observation above Lemma \ref{L:one}, $r=0$ or $r=d$.  Replacing $\{\theta_{i}\}_{i=0}^{d}$ with $\{\theta_{d-i}\}_{i=0}^{d}$ as necessary, we may assume, without loss of generality, that $r=0$.

Now $\theta_{0}$ is adjacent to an element of $\K$ not among $\{\theta_{i}\}_{i=0}^{d}$.  Recall this number is either $2-\theta_{0}$ or $-2-\theta_{0}$.  When $d\geq1$, let $\varepsilon$ be as below line (\ref{E:eigvals}).  For notational convenience we define $\varepsilon$ for $d=0$.  In this case if $\theta_{0}=\pm1$ we define $\varepsilon=\theta_{0}$ and if $\theta_{0}\ne\pm1$ we define $\varepsilon=1$.  For all values of $d$, $-2\varepsilon-\theta_{0}$ is not among $\{\theta_{i}\}_{i=0}^{d}$.  Therefore $V(-2\varepsilon-\theta_{0})=0$.  By this and Lemma \ref{L:tridi} we have $(y-\varepsilon z).V(\theta_{0})=0$.

We have that $\theta_{i}$ satisfies (\ref{E:eigvals}) for $0\leq i\leq d$.  For notational convenience, we define $\theta_{i}$ by the equation (\ref{E:eigvals}) for all integers $i\geq0$.

Let $0\ne w_{0}\in V(\theta_{0})$.  We have $(y-\varepsilon z).w_{0}=0$.  We define vectors \{$w_{i}\}_{i\geq1}$ recursively by
\begin{equation}\label{E:vi}
w_{i}=\frac{\varepsilon}{2i}(y+(-1)^{i-1}\varepsilon z).w_{i-1}\qquad i\geq 1.
\end{equation}
By Lemma \ref{L:tridi}, $w_{i}\in V(\theta_{i})$ for $i\geq 0$.  In (\ref{E:vi}) we replace $i$ with $i+1$ and rearrange the terms to get
\begin{equation}\label{E:viplus}
(y+(-1)^{i}\varepsilon z).w_{i}=2\varepsilon(i+1)w_{i+1}\qquad i\geq 0.
\end{equation}
We claim that, for $i\geq0$,
\begin{equation}\label{E:viminus}
(y-(-1)^{i}\varepsilon z).w_{i}=2\varepsilon(\varepsilon\theta_{0}-i+1)w_{i-1},
\end{equation}
where $w_{-1}=0$.  We do this using induction on $i$.  First assume $i=0$.  Then (\ref{E:viminus}) holds since both sides are equal to $0$.  Now assume $i\geq1$.  Using (\ref{E:rel2}) we check that $(y+z)^{2}-(y-z)^{2}=4x$.  This implies $(y+\varepsilon z)^{2}-(y-\varepsilon z)^{2}=\varepsilon4x$, so
\begin{equation}\label{E:squareval}
((y+\varepsilon z)^{2}-(y-\varepsilon z)^{2}).w_{i-1}=\varepsilon4x.w_{i-1}.
\end{equation}
As we evaluate (\ref{E:squareval}), we consider two cases:

Case 1 ($i$ is even): By (\ref{E:viplus}), we have $(y-\varepsilon z).w_{i-1}=2\varepsilon iw_{i}$ and $(y+\varepsilon z).w_{i-2}=2\varepsilon(i-1)w_{i-1}$.  By (\ref{E:viminus}) and induction we have $(y+\varepsilon z).w_{i-1}=2\varepsilon(\varepsilon\theta_{0}-i+2)w_{i-2}$.  By (\ref{E:eigvals}) we have $x.w_{i-1}=(2\varepsilon(i-1)-\theta_{0})w_{i-1}$.  By these comments and (\ref{E:squareval}) we routinely obtain (\ref{E:viminus}).

Case 2 ($i$ is odd): By (\ref{E:viplus}), we have $(y+\varepsilon z).w_{i-1}=2\varepsilon iw_{i}$ and $(y-\varepsilon z).w_{i-2}=2\varepsilon(i-1)w_{i-1}$.  By (\ref{E:viminus}) and induction we have $(y-\varepsilon z).w_{i-1}=2\varepsilon(\varepsilon\theta_{0}-i+2)w_{i-2}$.  By (\ref{E:eigvals}) we have $x.w_{i-1}=(\theta_{0}-2\varepsilon(i-1))w_{i-1}$.  By these comments and (\ref{E:squareval}) we routinely obtain (\ref{E:viminus}).

We have now verified (\ref{E:viminus}).  We next claim that, for $i\geq0$,
\begin{align}
x.w_{i}=&(-1)^{i}(\theta_{0}-2\varepsilon i)w_{i},\label{E:Gen1}\\
y.w_{i}=&\varepsilon(\varepsilon\theta_{0}-i+1)w_{i-1}+\varepsilon(i+1)w_{i+1},\label{E:Gen2}\\
z.w_{i}=&(-1)^{i-1}(\varepsilon\theta_{0}-i+1)w_{i-1}+(-1)^{i}(i+1)w_{i+1}.\label{E:Gen3}
\end{align}

Adding (\ref{E:viplus}), (\ref{E:viminus}) and dividing the result by $2$ we get (\ref{E:Gen2}).  Adding $(-1)^{i}\varepsilon$ times (\ref{E:viplus}) to $(-1)^{i-1}\varepsilon$ times (\ref{E:viminus}) and dividing the result by $2$ we get (\ref{E:Gen3}).  Combining the fact that $w_{i}\in V(\theta_{i})$ with (\ref{E:eigvals}) we obtain (\ref{E:Gen1}).  The claim follows.

By (\ref{E:Gen1})--(\ref{E:Gen3}), span$\{w_{i}\}_{i\geq0}$ is closed under the actions of $x,y,z$, and hence all of $\mathcal{A}$.  Because the $\mathcal{A}$-module $V$ is irreducible and $w_{0}\ne0$, we have span$\{w_{i}\}_{i\geq0}=V$.

We now show there exists a nonnegative integer $t$ such that $w_{t}=0$.  By construction, the sequences $\{\theta_{2i}\}_{i\geq0}$ and $\{\theta_{2i+1}\}_{i\geq0}$ are arithmetic progressions, so $\{\theta_{i}\}_{i\geq0}$ has an infinite number of distinct elements.  Since $V$ is finite-dimensional, there must be a nonnegative integer $i$ such that $\theta_{i}$ is not among $\{\theta_{j}\}_{j=0}^{d}$.  Observe $V(\theta_{i})=0$ so $w_{i}=0$.

Assume $w_{t}=0$.  We now show that $t\geq d+1$.  Assume $t\leq d$.  By (\ref{E:viplus}) $w_{i}=0$ for all $i\geq t$.  Therefore $V=\mathrm{span}\{w_{i}\}_{i=0}^{t-1}$.  By this, and the fact that $\{\theta_{i}\}_{i=0}^{d}$ are distinct, we have that $V(\theta_{d})=0$, a contradiction.  Therefore $t\geq d+1$.

Let $c$ denote the smallest integer such that $c\geq d$ and $w_{c+1}=0$.  Then $V=\mathrm{span}\{w_{i}\}_{i=0}^{c}$.  Setting $i=c+1$ in (\ref{E:viminus}) and using $w_{c+1}=0$, we get $2\varepsilon(\varepsilon\theta_{0}-c)w_{c}=0$, but $2\varepsilon$ and $w_{c}$ are nonzero, so $\varepsilon\theta_{0}-c=0$.  This means $\theta_{0}=\varepsilon c$.  In particular $\theta_{0}$ is an integer, so, by (\ref{E:eigvals}), $\theta_{i}$ are integers for all $i\geq0$.

Also by (\ref{E:eigvals}), $\{\theta_{i}\}_{i\geq0}$ are either all even or all odd.  We now consider these two subcases separately.

Case 1 ($\{\theta_{i}\}_{i\geq0}$ are even): Since $\theta_{d}$ is even, it is not equal to $\pm1$.  By rows $2$ and $4$ of the second table from Lemma \ref{L:one}, $\theta_{d}$ is adjacent to an element of $\K$ not among $\{\theta_{i}\}_{i=0}^{d}$.  Therefore $\theta_{d+1}$ is not among $\{\theta_{i}\}_{i=0}^{d}$.  This means $w_{d+1}=0$, so $c\leq d$.  By this and the fact that $c\geq d$, we have $c=d$.

From this we draw two conclusions.  First of all, using $\theta_{0}=\varepsilon c$, we find $\theta_{0}=\varepsilon d$.  Secondly, we find that the vectors  $\{w_{i}\}_{i=0}^{d}$ form a basis for $V$.  If $\varepsilon=1$, we define $v_{i}=w_{i}$ for $0\leq i\leq d$.  If $\varepsilon=-1$, we define $v_{i}=(-1)^{i}w_{d-i}$.  In both cases $\{v_{i}\}_{i=0}^{d}$ is a basis for $V$.  Combining the construction of $\{v_{i}\}_{i=0}^{d}$, (\ref{E:Gen1})--(\ref{E:Gen3}) and $\theta_{0}=\varepsilon d$, we obtain (\ref{E:Bd1})--(\ref{E:Bd3}).

Case 2 ($\{\theta_{i}\}_{i\geq0}$ are odd): Recall $\theta_{0}=\varepsilon c$ so $c$ is odd.  Therefore there exists an integer $k\geq0$ such that $c=2k+1$.  We show that $k=d$.  By (\ref{E:eigvals}) and since $c$ is odd we have $\theta_{i}=\theta_{c-i}$ for $0\leq i\leq c$.  In this equation we set $i=k$ to get $\theta_{k}=\theta_{k+1}$.  Because $\{\theta_{i}\}_{i=0}^{d}$ are distinct, $k\geq d$.

This implies $c\geq2d+1>d$, so $V(\theta_{d+1})\ne0$.  By Lemma \ref{L:one} rows $5$--$8$, we have $\theta_{d}=\pm1$.  By (\ref{E:eigvals}) with $i=d$, we get $\theta_{d}=(-1)^{d}(\varepsilon c-2\varepsilon d)$, so $c-2d=\pm1$.  This means $k$ is either $d$ or $d-1$, but $k\geq d$.  Therefore $k=d$ and hence $c=2d+1$, so $\theta_{i}=\varepsilon(-1)^{i}(2d-2i+1)$.

We now have $V=\mathrm{span}\{w_{i}\}_{i=0}^{2d+1}$.  Let $V_{0}=\mathrm{span}\{w_{i}+w_{c-i}\}_{i=0}^{d}$ and $V_{1}=\mathrm{span}\{w_{i}-w_{c-i}\}_{i=0}^{d}$.  Observe by construction that $V=V_{0}+V_{1}$ and by (\ref{E:Gen1})--(\ref{E:Gen3}), $V_{0}$ and $V_{1}$ are closed under the action of $\mathcal{A}$.  By these comments and the fact that the $\mathcal{A}$-module $V$ is irreducible and the fact that $V=V_{0}+V_{1}$, either $V_{0}=0$ and $V_{1}=V$, or $V_{1}=0$ and $V_{0}=V$.  In the former case, we define $\delta=-1$ and in the latter case, we define $\delta=1$.  Then $w_{i}=\delta w_{c-i}$ for $0\leq i\leq c$ and the vectors $\{w_{i}\}_{i=0}^{d}$ form a basis for $V$.  Let $v_{i}=\delta^{i} w_{i}$ for $0\leq i\leq c$.  Then $\{v_{i}\}_{i=0}^{d}$ is a basis for $V$ and $v_{i}=v_{c-i}$ for $0\leq i\leq c$.

Using (\ref{E:Gen1})--(\ref{E:Gen3}) and the definition of $\{v_{i}\}_{i=0}^{d}$, we determine the actions of $x,y,z$ on $\{v_{i}\}_{i=0}^{d}$ for the different values of $\varepsilon,\delta$.  Comparing these actions with the data from Lemmas \ref{L:ABd0}--\ref{L:ABdz}, we find that the $\mathcal{A}$-module $V$ is in the isomorphism class displayed in the table below.
\begin{center}
\begin{tabular}{|l||c|c|}
\hline
& $(-1)^{d}\delta=1$ & $(-1)^{d}\delta=-1$ \\ \hline\hline
$(-1)^{d}\varepsilon=1$ & $AB(d,0)$ & $AB(d,x)$ \\ \hline
$(-1)^{d}\varepsilon=-1$ & $AB(d,y)$ & $AB(d,z)$ \\ \hline
\end{tabular}
\end{center}
\hfill $\Box$ \\

We comment on Theorem \ref{T:class}.  In \cite{arik-kayser}, Arik and Kayserilioglu introduced a complex unital associative algebra with generators $J_{1},J_{2},J_{3}$ and relations
\begin{equation}\label{E:ACSA}
\{J_{1},J_{2}\}=J_{3},\qquad\{J_{2},J_{3}\}=J_{1},\qquad\{J_{3},J_{1}\}=J_{2},
\end{equation}
where $\{A,B\}=AB+BA$.  They called their algebra the \emph{anticommutator spin algebra}, abbreviated $ACSA$.  Comparing equations (\ref{E:rel1})--(\ref{E:rel3}) and (\ref{E:ACSA}), we see that, when $\K=\C$, there is an algebra isomorphism $\mathcal{A}\to ACSA$ that sends $x\mapsto2J_{3},y\mapsto2J_{1},z\mapsto2J_{2}$.  Arik and Kayserilioglu claimed to classify up to isomorphism all finite-dimensional irreducible representations of $ACSA$.  However, their result is incorrect; they only found three types of representations instead of the five described in Lemmas \ref{L:Bd}--\ref{L:ABdz}.  What Arik and Kayserilioglu actually classified were the possible eigenvalue sequences for the action of $J_{3}$ in a finite-dimensional irreducible representation.  But the distinct isomorphism classes $AB(d,0)$ and $AB(d,x)$ yield the same eigenvalue sequence for the action of $J_{3}$.  Similarly, the distinct isomorphism classes $AB(d,y)$ and $AB(d,z)$ yield the same eigenvalue sequence for the action of $J_{3}$.

We include a result for later use.

\begin{lemma}\label{L:htrace}
Let $V$ denote an $\mathcal{A}$-module contained in one of the five families from Lemmas \ref{L:Bd}--\ref{L:ABdz}.  Then for $n\in\I$, the trace of $h_{n}$ on $V$ is given on the following table.

\begin{center}
\begin{tabular}{|l||c|c|c|c|}
\hline
& $\mbox{\em{tr}}(h_{0})$ & $\mbox{\em{tr}}(h_{x})$ & $\mbox{\em{tr}}(h_{y})$ & $\mbox{\em{tr}}(h_{z})$ \\ \hline\hline
$B(d)$ & $0$ & $0$ & $0$ & $0$ \\ \hline
$AB(d,0)$ & $3(-1)^{d}(d+1)$ & $(-1)^{d+1}(d+1)$ & $(-1)^{d+1}(d+1)$ & $(-1)^{d+1}(d+1)$ \\ \hline
$AB(d,x)$ & $(-1)^{d+1}(d+1)$ & $3(-1)^{d}(d+1)$ & $(-1)^{d+1}(d+1)$ & $(-1)^{d+1}(d+1)$ \\ \hline
$AB(d,y)$ & $(-1)^{d+1}(d+1)$ & $(-1)^{d+1}(d+1)$ & $3(-1)^{d}(d+1)$ & $(-1)^{d+1}(d+1)$ \\ \hline
$AB(d,z)$ & $(-1)^{d+1}(d+1)$ & $(-1)^{d+1}(d+1)$ & $(-1)^{d+1}(d+1)$ & $3(-1)^{d}(d+1)$ \\ \hline
\end{tabular}
\end{center}
\end{lemma}

\noindent {\it Proof:}
Apply Theorem \ref{T:trace} to Definition \ref{D:h0123}.
\hfill $\Box$ \\

\section{The $G$-action on the $\mathcal{A}$-modules}\label{S:S4}

Recall the subgroup $G\subseteq$Aut$(\mathcal{A})$ from below Lemma \ref{L:KK4}.  Let $V$ denote a finite-dimensional irreducible $\mathcal{A}$-module.  In this section we show what happens when we twist $V$ via an element of $G$.

\begin{definition}\label{D:twisted}
{\normalfont Let $V$ denote an $\mathcal{A}$-module.  For $\sigma\in$Aut$(\mathcal{A})$ there exists an $\mathcal{A}$-module structure on $V$, called $V$ \emph{twisted via} $\sigma$ that behaves as follows: for all $a\in\mathcal{A},v\in V$, the vector $a.v$ computed in $V$ twisted via $\sigma$ coincides with the vector $\sigma^{-1}(a).v$ computed in the original $\mathcal{A}$-module $V$.  Sometimes we abbreviate ${}^{\sigma}V$ for $V$ twisted via $\sigma$.  Observe that Aut$(\mathcal{A})$ acts on the set of $\mathcal{A}$-modules, with $\sigma$ sending $V$ to ${}^{\sigma}V$ for all $\sigma\in$Aut$(\mathcal{A})$ and every $\mathcal{A}$-module $V$.  Observe that $V$ and ${}^{\sigma}V$ have the same dimension and that ${}^{\sigma}V$ is irreducible if and only if $V$ is irreducible.}
\end{definition}

In Section \ref{S:modules} we described the set of isomorphism classes of finite-dimensional irreducible $\mathcal{A}$-modules.  From Definition \ref{D:twisted}, $G$ acts on this set.  We now investigate this $G$-action.  Recall from Definition \ref{D:I} that $\I$ consists of the symbols $0,x,y,z$.

\begin{theorem}\label{T:modact}
Let $V$ denote a finite-dimensional irreducible $\mathcal{A}$-module of diameter $d$ and let $\sigma\in G$.  Then the following (i), (ii) hold.
\begin{enumerate}
\item[\rm (i)] Assume $V$ is of type $B(d)$.  Then ${}^{\sigma}V$ is of type $B(d)$.

\item[\rm (ii)] Assume $V$ is of type $AB(d,n)$ for some $n\in\I$.  Then ${}^{\sigma}V$ is of type $AB(d,\tilde{\sigma}(n))$,
\end{enumerate}
where $\tilde{\sigma}$ is from Proposition \ref{P:S4action}.
\end{theorem}

\noindent {\it Proof:}
For $m\in\I$, the action of $h_{m}$ on ${}^{\sigma}V$ coincides with the action of $\sigma^{-1}(h_{m})$ on $V$.  Therefore, the trace of $h_{m}$ on ${}^{\sigma}V$ is equal to the trace of $\sigma^{-1}(h_{m})$ on $V$.  We evaluate the table from Lemma \ref{L:htrace} using this and Proposition \ref{P:S4action}.  The result follows.
\hfill $\Box$ \\

By Theorem \ref{T:modact}(i) the isomorphism class $B(d)$ is stabilized by everything in $G$.  For $n\in\I$ we now determine the stabilizer in $G$ of the isomorphism class of type $AB(d,n)$.  Recall the subgroups $K,S\subseteq G$ from below Lemma \ref{L:KK4}.

\begin{definition}\label{D:rhos}
{\normalfont Recall that the group $K$ consists of the automorphisms of $\mathcal{A}$ that fix each of $x,y,z$ up to sign.  Recall that $|K|=4$ by Lemma \ref{L:KK4}.  We define a bijection $\I\to K,n\mapsto\rho_{n}$ as follows.  The automorphism $\rho_{0}$ is the identity element of $K$.  For each nonzero $n\in\I$, by Lemma \ref{L:KK4}, there exists a unique element of $K$ that fixes $n$ and changes the sign of the other two elements of $\{x,y,z\}$.  We denote this element of $K$ by $\rho_{n}$.}
\end{definition}

Recall the group $\tilde{G}$ of permutations of $\I$ and the isomorphism $G\to\tilde{G}$ from Proposition \ref{P:S4action}.  Note that, for nonzero $n\in\I$, $\tilde{\rho}_{n}=(0n)(ml)$ where $m,l$ are the remaining nonzero elements of $\I$.

\begin{lemma}\label{L:stabilizers}
Let $n\in\I$ and let $V$ denote a finite-dimensional irreducible $\mathcal{A}$-module of type $AB(d,n)$.  Then, for $\sigma\in G$, the following (i)--(iii) are equivalent.
\begin{enumerate}
\item[\rm (i)] ${}^{\sigma}V$ is of type $AB(d,n)$.

\item[\rm (ii)] $\tilde{\sigma}$ fixes $n$.

\item[\rm (iii)] $\sigma\in\rho_{n}S\rho_{n}^{-1}$, where $\rho_{n}$ is from Definition \ref{D:rhos}.
\end{enumerate}
\end{lemma}

\noindent {\it Proof:}
(i)$\Leftrightarrow$(ii): Follows from Proposition \ref{P:S4action} and Theorem \ref{T:modact}.

(ii)$\Leftrightarrow$(iii): First assume that $n=0$, so that $\rho_{n}$ is the identity.  Then $\rho_{n}S\rho_{n}^{-1}=S$.  By Lemma \ref{L:S3}, $\tilde{S}$ consists of the permutations of $\I$ that fix $0$.  Now assume $n\ne0$.  Then, by Lemma \ref{L:S3} and the note after Definition \ref{D:rhos}, we check that $\tilde{\rho}_{n}\tilde{S}\tilde{\rho}^{-1}_{n}$ consists of the permutations of $\I$ that fix $n$.  Therefore $\sigma\in\rho_{n}S\rho^{-1}_{n}$ if and only if $\tilde{\sigma}$ fixes $n$.
\hfill $\Box$ \\

\section{The primitive idempotents}\label{S:idempotents}

In this section we determine the eigenvalues for the actions of $x,y,z$ on a finite-dimensional irreducible $\mathcal{A}$-module, and we define the corresponding primitive idempotents.

\begin{definition}\label{D:e}
{\normalfont Let $V$ denote a vector space over $\K$ with positive finite dimension and let $b:V\to V$ denote a diagonalizable linear transformation.  Let  $\{V_{i}\}_{i=0}^{d}$ denote an ordering of the eigenspaces of $b$.  For $0\leq i\leq d$ let $\theta_{i}$ denote the eigenvalue for $b$ associated with $V_{i}$ and define $e_{i}\in\rm{End}(V)$ such that $(e_{i}-I)V_{i}=0$ and $e_{i}V_{j}=0$ for $j\ne i$ $(0\leq j\leq d)$.  Here $I$ denotes the identity of End$(V)$.  We call $e_{i}$ the \emph{primitive idempotent of $b$} corresponding to $\theta_{i}$.  Observe that
\begin{itemize}
\item[\rm (i)] $\sum_{i=0}^{d}e_{i}=I$,

\item[\rm (ii)] $e_{i}e_{j}=\delta_{ij}e_{i}$\qquad$(0\leq i,j\leq d)$,

\item[\rm (iii)] $ae_{i}=\theta_{i}e_{i}$\qquad$(0\leq i\leq d)$,

\item[\rm (iv)] $e_{i}V=V_{i}$\qquad$(0\leq i\leq d)$.
\end{itemize}}
\end{definition}

Note that
\begin{equation}\label{E:e}
e_{i}=\prod_{
\genfrac{}{}{0pt}{} {0\leq j\leq d}{j\ne i}
}\frac{b-\theta_{j}I}{\theta_{i}-\theta_{j}}\qquad(0\leq i\leq d).
\end{equation}

We will now determine the eigenvalues for the actions of $x,y,z$ on a finite-dimensional irreducible $\mathcal{A}$-module $V$.  To do this, we will first determine the eigenvalues for the action of $x,y,z$ when $V$ is of type $B(d)$ or $AB(d,0)$.  Then we will determine the eigenvalues for the actions of $x,y,z$ when $V$ is of type $AB(d,n)$ for nonzero $n\in\I$.

\begin{proposition}\label{P:idemtable1}
Let $V$ denote a finite-dimensional irreducible $\mathcal{A}$-module of type $B(d)$ or $AB(d,0)$.  For each of $x,y,z$ the action on $V$ is diagonalizable.  The eigenvalues for this action are given in the table below.
\begin{center}
\begin{tabular}{|l||c|c|}
\hline
& $B(d)$ & $AB(d,0)$ \\ \hline\hline
$x$ & $(-1)^{i}(d-2i)$ & $(-1)^{d+i}(2d-2i+1)$ \\ \hline
$y$ & $(-1)^{i}(d-2i)$ & $(-1)^{d+i}(2d-2i+1)$ \\ \hline
$z$ & $(-1)^{i}(d-2i)$ & $(-1)^{d+i}(2d-2i+1)$ \\ \hline
\end{tabular}
\end{center}
In the above table, the integer $i$ runs from $0$ to $d$.
\end{proposition}

\noindent {\it Proof:} 
If $V$ is of type $B(d)$, then by Lemma \ref{L:Bd}, the action of $x$ on $V$ is diagonalizable with the desired eigenvalues.  If $V$ is of type $AB(d,0)$, then by Lemma \ref{L:ABd0} the action of $x$ on $V$ is diagonalizable with the desired eigenvalues.  We have now verified our assertions for $x$.

We now verify our assertions for $y,z$.  To that end, let $a$ denote one of $y,z$.  Pick an element $\sigma\in S$ such that $\sigma(a)=x$.  By Theorem \ref{T:modact} and Lemma \ref{L:stabilizers}, the twisted module ${}^{\sigma}V$ is of the same type as $V$.  By Definition \ref{D:twisted}, the action of $x$ on ${}^{\sigma}V$ coincides with the action of $\sigma^{-1}(x)=a$ on the untwisted module $V$.  Therefore the actions are both diagonalizable and have the same eigenvalues.
\hfill $\Box$ \\

\begin{proposition}\label{P:idemtable2}
Fix a nonzero $n\in\I$ and let $V$ denote a finite-dimensional irreducible $\mathcal{A}$-module of type $AB(d,n)$.  For each of $x,y,z$ the action on $V$ is diagonalizable.  The eigenvalues for this action are given in the table below.
\begin{center}
\begin{tabular}{|l||c|c|c|c|c|}
\hline
& $AB(d,x)$ & $AB(d,y)$ & $AB(d,z)$ \\ \hline\hline
$x$ & $(-1)^{d+i}(2d-2i+1)$ & $(-1)^{d+i+1}(2d-2i+1)$ & $(-1)^{d+i+1}(2d-2i+1)$ \\ \hline
$y$ & $(-1)^{d+i+1}(2d-2i+1)$ & $(-1)^{d+i}(2d-2i+1)$ & $(-1)^{d+i+1}(2d-2i+1)$ \\ \hline
$z$ & $(-1)^{d+i+1}(2d-2i+1)$ & $(-1)^{d+i+1}(2d-2i+1)$ & $(-1)^{d+i}(2d-2i+1)$ \\ \hline
\end{tabular}
\end{center}
In the above table, the integer $i$ runs from $0$ to $d$.
\end{proposition}

\noindent {\it Proof:} 
Recall the automorphism $\rho=\rho_{n}$ of $\mathcal{A}$ from Definition \ref{D:rhos}.  By Theorem \ref{T:modact} and the note at the end of Definition \ref{D:rhos} we find that the twisted module ${}^{\rho}V$ is of type $AB(d,0)$.  Let $a$ denote one of $x,y,z$ and note that, by Proposition \ref{P:idemtable1}, the action of $a$ on ${}^{\rho}V$ is diagonalizable with eigenvalues $\{(-1)^{d+i}(2d-2i+1)\}_{i=0}^{d}$.  By Definition \ref{D:twisted}, the action of $\rho(a)$ on the untwisted module $V$ is diagonalizable with eigenvalues $\{(-1)^{d+i}(2d-2i+1)\}_{i=0}^{d}$.  By Definition \ref{D:rhos}, $\rho(a)=a$ when $n=a$ and $\rho(a)=-a$ when $n\ne a$.  The result follows.
\hfill $\Box$ \\

\begin{definition}\label{D:eieio1}
{\normalfont Let $V$ denote a finite-dimensional irreducible $\mathcal{A}$-module of type $B(d)$ or $AB(d,0)$.  For $a$ among $x,y,z$ and $0\leq i\leq d$, let $\theta_{i}^{a}$ denote the $i^{\mathrm{th}}$ eigenvalue of $a$ on $V$ from the table in Proposition \ref{P:idemtable1}.  We define $e_{i}^{a}$ to be the primitive idempotent associated with $\theta_{i}^{a}$, for the action of $a$ on $V$.}
\end{definition}

\begin{definition}\label{D:eieio2}
{\normalfont For nonzero $n\in\I$, let $V$ denote a finite-dimensional irreducible $\mathcal{A}$-module of type $AB(d,n)$.  For $a$ among $x,y,z$ and $0\leq i\leq d$, let $\theta_{i}^{a}$ denote the $i^{\mathrm{th}}$ eigenvalue of $a$ on $V$ from the table in Proposition \ref{P:idemtable2}.  We define $e_{i}^{a}$ to be the primitive idempotent associated with $\theta_{i}^{a}$, for the action of $a$ on $V$.}
\end{definition}

Recall the notion of standard order from below Corollary \ref{C:path}.

\begin{lemma}\label{L:stord}
Let $a$ be among $x,y,z$.  With respect to Definitions \ref{D:eieio1}, \ref{D:eieio2}, the ordering $\{\theta_{i}^{a}\}_{i=0}^{d}$ is standard.
\end{lemma}

\noindent {\it Proof:} 
Use the tables in Propositions \ref{P:idemtable1}, \ref{P:idemtable2}.
\hfill $\Box$ \\

We now present two slightly technical results that will be used in later sections.

\begin{lemma}\label{L:Bsigma}
Let $V$ denote a finite-dimensional irreducible $\mathcal{A}$-module of type $B(d)$.  Pick an element $\sigma\in S$.  Pick $a$ among $x,y,z$.  Then $\sigma(e_{i}^{a})=e_{i}^{\sigma(a)}$ for $0\leq i\leq d$. 
\end{lemma}

\noindent {\it Proof:} 
The idempotents $e_{i}^{a},e_{i}^{\sigma(a)}$ are found using (\ref{E:e}).  By Proposition \ref{P:idemtable1}, $\theta_{i}^{a}=\theta_{i}^{\sigma(a)}$.  The result follows.
\hfill $\Box$ \\

We set some notation for later use.  Let $0\ne a\in\I$.  We define the function $\widehat{a}:\I\to\K$ to by $\widehat{a}(n)=1$ for $n\in\{0,a\}$ and $\widehat{a}(n)=-1$ for $n\in\I\setminus\{0,a\}$.

\begin{lemma}\label{L:ABsigma}
Let $V$ denote a finite-dimensional irreducible $\mathcal{A}$-module of type $AB(d,n)$ and let $\rho_{n}$ be as in Definition \ref{D:rhos}.  Pick an element $\sigma\in G$ and let $\tau=\rho_{n}\sigma\rho_{n}^{-1}$.  Pick $a$ among $x,y,z$.  Then $\tau(e_{i}^{a})=e_{i}^{\sigma(a)}$ for $0\leq i\leq d$. 
\end{lemma}

\noindent {\it Proof:} 
The idempotents $e_{i}^{a},e_{i}^{\sigma(a)}$ are found using (\ref{E:e}).  We have $\tau(a)=\widehat{a}(n)\widehat{\sigma(a)}(n)\sigma(a)$ by Definition \ref{D:rhos} and $\theta_{i}^{a}=\widehat{a}(n)\widehat{\sigma(a)}(n)\theta_{i}^{\sigma(a)}$ by Propositions \ref{P:idemtable1}, \ref{P:idemtable2}.  The result follows.
\hfill $\Box$ \\

\section{Six bases for $V$}\label{S:bases}

Let $V$ denote a finite-dimensional irreducible $\mathcal{A}$-module.  In this section we will display six bases for $V$ with respect to which the matrices representing $x,y,z$ are attractive.  To begin, we will look at the basis for $V$ provided in Lemmas \ref{L:Bd}--\ref{L:ABdz}.

\begin{lemma}\label{L:xybasis}
Let $V$ denote a finite-dimensional irreducible $\mathcal{A}$-module.  Let $\{v_{i}\}_{i=0}^{d}$ denote the basis for $V$ from Lemmas \ref{L:Bd}--\ref{L:ABdz}.  Then the following (i)--(iii) hold.
\begin{itemize}
\item[\rm (i)] $v_{i}\in e_{i}^{x}V\qquad(0\leq i\leq d)$.

\item[\rm (ii)] Let $v=\sum_{i=0}^{d}v_{i}$.  Then $v\in e_{0}^{y}V$.

\item[\rm (iii)] $v_{i}=e_{i}^{x}v\qquad(0\leq i\leq d)$.
\end{itemize}
\end{lemma}

\noindent {\it Proof:} 
(i) Follows from equations (\ref{E:Bd1}), (\ref{E:ABdx1}), (\ref{E:ABdy1}), (\ref{E:ABdz1}) and Propositions \ref{P:idemtable1}, \ref{P:idemtable2}.

(ii) Follows from equations (\ref{E:Bd2}), (\ref{E:ABdx2}), (\ref{E:ABdy2}), (\ref{E:ABdz2}) and Propositions \ref{P:idemtable1}, \ref{P:idemtable2}.

(iii) By part (ii), $e_{i}v=\sum_{j=0}^{d}e_{i}^{x}v_{j}$.  Since $v_{j}\in e_{j}^{x}V$ we have $e_{i}^{x}v_{j}=\delta_{ij}v_{i}$.  The result follows.
\hfill $\Box$ \\

\begin{lemma}\label{L:eie0}
Let $V$ denote a finite-dimensional irreducible $\mathcal{A}$-module.  Pick $a,b$ among $x,y,z$ with $a\ne b$.  Then the action of $e_{i}^{a}e_{0}^{b}$ on $V$ is nonzero for $0\leq i\leq d$.
\end{lemma}

\noindent {\it Proof:} 
Observe $e_{i}^{x}e_{0}^{y}V$ is nonzero because it contains the nonzero vector $v_{i}$ from Lemma \ref{L:xybasis}.  Now, let $\sigma\in S$ denote the unique automorphism of $\mathcal{A}$ such that $\sigma(a)=x$ and $\sigma(b)=y$.  Let $\rho\in K$ denote the identity if $V$ is of type $B(d)$ and $\rho_{n}$ if $V$ is of type $AB(d,n)$.  Let $\tau=\rho\sigma\rho^{-1}$.  By Lemma \ref{L:stabilizers}, the $\mathcal{A}$-modules $V$ and ${}^{\tau}V$ are isomorphic, so the action of $e_{i}^{x}e_{0}^{y}$ on ${}^{\tau}V$ is nonzero.  By Lemmas \ref{L:Bsigma}, \ref{L:ABsigma}, the action of $e_{i}^{x}e_{0}^{y}$ on ${}^{\tau}V$ coincides with the action of $e_{i}^{a}e_{0}^{b}$ on $V$.  Therefore $e_{i}^{a}e_{0}^{b}V\ne0$ as desired.
\hfill $\Box$ \\

We now obtain six bases for $V$.

\begin{theorem}\label{T:bases}
Let $V$ denote a finite-dimensional irreducible $\mathcal{A}$-module of diameter $d$.  Pick $a,b$ among $x,y,z$ with $a\ne b$.  Then, for $0\ne v^{b}\in e_{0}^{b}V$ and $0\leq i\leq d$, $e_{i}^{a}v^{b}$ is nonzero and therefore a basis for $e_{i}^{a}V$.  Moreover, the sequence $\{e_{i}^{a}v^{b}\}_{i=0}^{d}$ is a basis for $V$.
\end{theorem}

\noindent {\it Proof:} 
We have dim$(e_{0}^{b}V)=1$ and $0\ne v^{b}\in e_{0}^{b}V$, so $v^{b}$ spans $e_{0}^{b}V$.  Therefore $e_{i}^{a}v^{b}$ spans $e_{i}^{a}e_{0}^{b}V$.  Now $e_{i}^{a}v^{b}\ne0$ in view of Lemma \ref{L:eie0}.
\hfill $\Box$ \\

\section{The matrices representing $x,y,z$ with respect to the six bases}\label{S:matrices}

Let $V$ denote a finite-dimensional irreducible $\mathcal{A}$-module.  In Theorem \ref{T:bases} we displayed six bases for $V$.  In this section we will display the matrices representing $x,y,z$ with respect to these bases.

\begin{lemma}\label{L:Babcrels}
Let $V$ denote a finite-dimensional irreducible $\mathcal{A}$-module of type $B(d)$.  Let $a,b,c$ denote a permutation of $x,y,z$.  For $0\leq i\leq d$, the following equations hold on $V$:
\begin{align}
ae_{i}^{a}e_{0}^{b}=&(-1)^{i}(d-2i)e_{i}^{a}e_{0}^{b},\label{E:Babc1}\\
be_{i}^{a}e_{0}^{b}=&(d-i+1)e_{i-1}^{a}e_{0}^{b}+(i+1)e_{i+1}^{a}e_{0}^{b},\label{E:Babc2}\\
ce_{i}^{a}e_{0}^{b}=&(-1)^{i-1}(d-i+1)e_{i-1}^{a}e_{0}^{b}+(-1)^{i}(i+1)e_{i+1}^{a}e_{0}^{b}.\label{E:Babc3}
\end{align} 
Here $e_{-1}^{a}=0$ and $e_{d+1}^{a}=0$.
\end{lemma}

\noindent {\it Proof:} 
For the case $(a,b,c)=(x,y,z)$ the equations (\ref{E:Babc1})--(\ref{E:Babc3}) are reformations of (\ref{E:Bd1})--(\ref{E:Bd3}) in light of Lemma \ref{L:xybasis}.  The remaining cases follow from Lemma \ref{L:Bsigma}.
\hfill $\Box$ \\

\begin{theorem}\label{T:Brep}
Let $V$ denote a finite-dimensional irreducible $\mathcal{A}$-module of type $B(d)$.  Pick $a,b$ among $x,y,z$ with $a\ne b$ and recall the basis $\{e_{i}^{a}v^{b}\}_{i=0}^{d}$ from Theorem \ref{T:bases}.  With respect to this basis, the matrices representing $x,y,z$ are described below.  Let $c$ denote the element of $\{x,y,z\}$ other than $a,b$.  The matrices are
\[
a:\mbox{\emph{diag}}(d,2-d,d-4,\ldots,4-d,d-2,-d),
\]

\[
b:
\begin{pmatrix}
0 & d & & & & & & \\
1 & 0 & d-1 & & & & \mathbf{0} &\\
& 2 & 0 & d-2 &  & & &\\
& & 3 & . & . &  & &\\
& & & . & . & . & &\\
& & & & . & . & 2 &\\
& \mathbf{0} & & & & d-1 & 0 & 1\\
& & & & & & d & 0
\end{pmatrix},
\]

\[
c:
\begin{pmatrix}
0 & d & & & & & & \\
1 & 0 & 1-d & & & & \mathbf{0} &\\
& -2 & 0 & d-2 & & & &\\
& & 3 & . & . & & &\\
& & & . & . & . & &\\
& & & & . & . & 2 &\\
& \mathbf{0} & & & & d-1 & 0 & -1\\
& & & & & & -d & 0
\end{pmatrix}.
\]
\end{theorem}

\noindent {\it Proof:} 
The actions for $a,b,c$ on $\{e_{i}^{a}v^{b}\}_{i=0}^{d}$ are found by applying equations (\ref{E:Babc1})--(\ref{E:Babc3}) to $v^{b}$ and recalling that $e_{0}^{b}v^{b}=v^{b}$.
\hfill $\Box$ \\

\begin{lemma}\label{L:ABabcrels}
Let $V$ denote a finite-dimensional irreducible $\mathcal{A}$-module of type $AB(d,n)$.  Let $a,b,c$ denote a permutation of $x,y,z$.  For $0\leq i\leq d$, the following equations hold on $V$:
\begin{align}
ae_{i}^{a}e_{0}^{b}=&\widehat{a}(n)(-1)^{d+i}(2d-2i+1)e_{i}^{a}e_{0}^{b},\label{E:ABabc1}\\
be_{i}^{a}e_{0}^{b}=&\widehat{b}(n)(-1)^{d}(2d-i+2)e_{i-1}^{a}e_{0}^{b}+\widehat{a}(n)(-1)^{d}(i+1)e_{i+1}^{a}e_{0}^{b},\label{E:ABabc2}\\
ce_{i}^{a}e_{0}^{b}=&\widehat{c}(n)(-1)^{i-1}(2d-i+2)e_{i-1}^{a}e_{0}^{b}+\widehat{c}(n)(-1)^{i}(i+1)e_{i+1}^{a}e_{0}^{b}.\label{E:ABabc3}
\end{align} 
Here $e_{-1}^{a}=0$ and $e_{d+1}^{a}=e_{d}^{a}$.  We are using the hat notation from above Lemma \ref{L:ABsigma}.
\end{lemma}

\noindent {\it Proof:} 
For the case $(a,b,c)=(x,y,z)$, the equations (\ref{E:ABabc1})--(\ref{E:ABabc3}) are reformations of (\ref{E:ABd01})--(\ref{E:ABdz3}) in light of Lemma \ref{L:xybasis}.  The remaining cases follow from Lemma \ref{L:ABsigma}.
\hfill $\Box$ \\

\begin{theorem}\label{T:ABrep}
Let $V$ denote a finite-dimensional irreducible $\mathcal{A}$-module of type $AB(d,n)$.  Pick $a,b$ among $x,y,z$ with $a\ne b$ and recall the basis $\{e_{i}^{a}v^{b}\}_{i=0}^{d}$ from Theorem \ref{T:bases}.  With respect to this basis, the matrices representing $x,y,z$ are described below.  Let $c$ denote the element of $\{x,y,z\}$ other than $a,b$.  The matrices are
\[
a:\widehat{a}(n)\mbox{\emph{diag}}((-1)^{d}(2d+1),\ldots,9,-7,5,-3,1),
\]

\[
b:
\widehat{b}(n)(-1)^{d}\begin{pmatrix}
0 & 2d+1 & & & & & & \\
1 & 0 & 2d & & & & \mathbf{0} &\\
& 2 & 0 & 2d-1 & & &\\
& & 3 & . & . & & &\\
& & & . & . & . & &\\
& & & & . & . & d+3 &\\
& \mathbf{0} & & & & d-1 & 0 & d+2\\
& & & & & & d & d+1
\end{pmatrix},
\]

\[
c:
\widehat{c}(n)\begin{pmatrix}
0 & 2d+1 & & & & & & \\
1 & 0 & -2d & & & & \mathbf{0} &\\
& -2 & 0 & 2d-1 & & & &\\
& & 3 & . & . & & &\\
& & & . & . & . & &\\
& & & & . & . & \mbox{\tiny $(-1)^{d-2}(d+3)$} &\\
& \mathbf{0} & & & & \mbox{\tiny $(-1)^{d-2}(d-1)$} & 0 & \mbox{\tiny $(-1)^{d-1}(d+2)$}\\
& & & & & & \mbox{\tiny $(-1)^{d-1}d$} & \mbox{\tiny $(-1)^{d}(d+1)$}
\end{pmatrix}.
\]
\end{theorem}

\noindent {\it Proof:} 
The actions for $a,b,c$ on $\{e_{i}^{a}v^{b}\}_{i=0}^{d}$ are found by applying equations (\ref{E:ABabc1})--(\ref{E:ABabc3}) to $v^{b}$ and recalling that $e_{0}^{b}v^{b}=v^{b}$.
\hfill $\Box$ \\

\begin{theorem}\label{T:Ltrip}
Let $V$ denote a finite-dimensional irreducible $\mathcal{A}$-module.  Then the actions of $x,y,z$ on $V$ form a Leonard triple.  If $V$ is of type $B$, then the Leonard triple is totally bipartite, and if $V$ is of type $AB$, then the Leonard triple is totally almost bipartite.
\end{theorem}

\noindent {\it Proof:} 
Use Definitions \ref{D:LT}, \ref{D:bipartiteLT} and the data from Theorems \ref{T:Brep}, \ref{T:ABrep}.
\hfill $\Box$ \\

\begin{corollary}\label{C:Ltrip}
Let $V$ denote a finite-dimensional irreducible $\mathcal{A}$-module.  For any nonzero scalars $\xi,\xi^{*},\xi^{\varepsilon}\in\K$, let $A=\xi x$, $A^{*}=\xi^{*}y$, $A^{\varepsilon}=\xi^{\varepsilon}z$  Then the actions of $A,A^{*},A^{\varepsilon}$ form a Leonard triple.  If $V$ is of type $B$, then the Leonard triple is totally bipartite, and if $V$ is of type $AB$, then the Leonard triple is totally almost bipartite.
\end{corollary}

\noindent {\it Proof:} 
Immediate.
\hfill $\Box$ \\

\section{Totally B/AB Leonard pairs of Bannai/Ito type}\label{S:LPclassification}

In Theorem \ref{T:Ltrip} and Corollary \ref{C:Ltrip} we displayed totally B/AB Leonard triples arising from finite-dimensional irreducible $\mathcal{A}$-modules.  In this section we classify the Leonard pairs associated with these Leonard triples.  We show that they correspond to a family of totally B/AB Leonard pairs said to have Bannai/Ito type.  Using this correspondence we classify the totally B/AB Leonard pairs of Bannai/Ito type with diameter at least $3$.

\begin{notation}\label{N:notation}
Let $V$ denote a vector space over $\K$ with finite positive dimension.  Let $A,A^{*}$ denote a Leonard pair on $V$.  Let $\{v_{i}\}_{i=0}^{d}$ denote a basis for $V$ with respect to which $A$ is diagonal and $A^{*}$ is irreducible tridiagonal. Let $\{v^{*}_{i}\}_{i=0}^{d}$ denote a basis for $V$ with respect to which $A^{*}$ is diagonal and $A$ is irreducible tridiagonal.  For $0\leq i\leq d$, let $\theta_{i}$ denote the eigenvalue for $A$ associated with $v_{i}$ and let $\theta_{i}^{*}$ denote the eigenvalue for $A^{*}$ associated with $v_{i}^{*}$.
\end{notation}

\begin{lemma}{\normalfont\cite[Theorem 1.5]{terw-vid}}\label{L:terw-vid}
With reference to Notation \ref{N:notation}, there exists a sequence of scalars $\beta,\gamma,\gamma^{*},\varrho,\varrho^{*},\omega,\eta,\eta^{*}$ taken from $\K$ such that both
\begin{align}
A^{2}A^{*}-\beta AA^{*}A+A^{*}A^{2}-\gamma(AA^{*}+A^{*}A)-\varrho A^{*}&=\gamma^{*}A^{2}+\omega A+\eta I,\label{E:AW1}\\
A^{*2}A-\beta A^{*}AA^{*}+AA^{*2}-\gamma^{*}(A^{*}A+AA^{*})-\varrho^{*}A&=\gamma A^{*2}+\omega A^{*}+\eta^{*} I.\label{E:AW2}
\end{align}
The sequence is uniquely determined by the pair $A,A^{*}$ provided the diameter is at least $3$.
\end{lemma}

The equations (\ref{E:AW1}), (\ref{E:AW2}) are known as the \emph{Askey-Wilson} relations \cite{Zhedanov}; see \cite{terw-vid}.

\begin{lemma}{\normalfont\cite[Theorem 1.9(v)]{terwilliger}}\label{L:formulae}
With reference to Notation \ref{N:notation}, the expressions
\begin{equation}\label{E:3trec}
\frac{\theta_{i-2}-\theta_{i+1}}{\theta_{i-1}-\theta_{i}},\qquad\qquad\frac{\theta^{*}_{i-2}-\theta^{*}_{i+1}}{\theta^{*}_{i-1}-\theta^{*}_{i}}
\end{equation}
are equal and independent of $i$ for $2\leq i\leq d-1$.
\end{lemma}

\begin{definition}{\normalfont\cite[Definition 4.2]{terw-vid}}\label{D:Plambda}
{\normalfont Given scalars $\beta,\gamma,\varrho$ in $\K$ we define a two-variable polynomial
\[
P(\lambda,\mu)=\lambda^{2}-\beta\lambda\mu+\mu^{2}-\gamma(\lambda+\mu)-\varrho.
\]
Given scalars $\beta,\gamma^{*},\varrho^{*}$ in $\K$ we define a two-variable polynomial
\[
P^{*}(\lambda,\mu)=\lambda^{2}-\beta\lambda\mu+\mu^{2}-\gamma^{*}(\lambda+\mu)-\varrho^{*}.
\]}
\end{definition}

We introduce further notation.

\begin{notation}\label{N:ais}
With reference to Notation \ref{N:notation}, for $0\leq i\leq d$, let $a_{i}$ (resp. $a_{i}^{*}$) denote the $(i,i)$-entries for the matrix representing $A$ (resp. $A^{*}$) with respect to the basis $\{v^{*}_{i}\}_{i=0}^{d}$ (resp. $\{v_{i}\}_{i=0}^{d}$).
\end{notation}

We obtain some formulae involving $\{a_{i}\}_{i=0}^{d},\{a_{i}^{*}\}_{i=0}^{d}$.

\begin{lemma}{\normalfont\cite[Corollary 5.2]{terw-vid}}\label{L:formulae2}
Let $\beta,\gamma,\gamma^{*},\varrho,\varrho^{*},\omega,\eta,\eta^{*}$ denote scalars in $\K$.  Then with reference to Notation \ref{N:notation}, \ref{N:ais} and Definition \ref{D:Plambda}, the following (i), (ii) are equivalent.
\begin{itemize}
\item[\rm (i)] The sequence $\beta,\gamma,\gamma^{*},\varrho,\varrho^{*},\omega,\eta,\eta^{*}$ satisfies (\ref{E:AW1}) and (\ref{E:AW2}).

\item[\rm (ii)] For $1\leq i\leq d$ both
\begin{equation}\label{E:TV1}
P(\theta_{i-1},\theta_{i})=0,\qquad P^{*}(\theta_{i-1}^{*},\theta_{i}^{*})=0,
\end{equation}
and for $0\leq i\leq d$ both
\begin{align}
a_{i}^{*}P(\theta_{i},\theta_{i})=&\gamma^{*}\theta_{i}^{2}+\omega\theta_{i}+\eta,\label{E:TV2}\\
a_{i}P^{*}(\theta^{*}_{i},\theta^{*}_{i})=&\gamma\theta_{i}^{*2}+\omega\theta^{*}_{i}+\eta^{*}.\label{E:TV3}
\end{align}
\end{itemize}
\end{lemma}

Let the Leonard pair $A,A^{*}$ be from Notation \ref{N:notation}.  Observe that $A,A^{*}$ is bipartite (resp. dual bipartite) if and only if $a_{i}$ (resp. $a_{i}^{*}$) is equal to $0$ for $0\leq i\leq d$.  Similarly, $A,A^{*}$ is almost bipartite (resp. dual almost bipartite) if and only if exactly one of $a_{0},a_{d}$ (resp. $a_{0}^{*},a_{d}^{*}$) is nonzero and $a_{i}$ (resp. $a_{i}^{*}$) is equal to $0$ for $1\leq i\leq d-1$.

\begin{lemma}\label{L:bipartite}
With reference to Notation \ref{N:notation}, the following (i), (ii) hold.  For $0\leq i\leq d$,
\begin{itemize}
\item[\rm (i)] Suppose $A,A^{*}$ is bipartite.  Then $\theta_{i}=-\theta_{d-i}$.

\item[\rm (ii)] Suppose $A,A^{*}$ is dual bipartite.  Then $\theta_{i}^{*}=-\theta_{d-i}^{*}$.
\end{itemize}
\end{lemma}

\noindent {\it Proof:} 
(i) Recall the bases $\{v_{i}\}_{i=0}^{d}$ and $\{v_{i}^{*}\}_{i=0}^{d}$ for $V$ from Notation \ref{N:notation}.  Let $s^{*}\in\rm{End}$$(V)$ be defined by $s^{*}.v^{*}_{i}=(-1)^{i}v^{*}_{i}$ for $0\leq i\leq d$.  By construction, $s^{*}$ is invertible, so $\{s^{*}.v_{i}\}_{i=0}^{d}$ is a basis for $V$.  Because the matrix representing $A$ with respect to the basis $\{v_{i}^{*}\}_{i=0}^{d}$ is bipartite tridiagonal, we have $As^{*}=-s^{*}A$.  Recall that $v_{i}$ is an eigenvector for $A$ with eigenvalue $\theta_{i}$ for $0\leq i\leq d$.  From these facts, we have that $s^{*}.v_{i}$ is an eigenvector for $A$ with eigenvalue $-\theta_{i}$ for $0\leq i\leq d$.  Therefore the matrix representing $A$ with respect to the basis $\{s^{*}.v_{i}\}_{i=0}^{d}$ is diagonal.  Because the matrix representing $A^{*}$ with respect to the basis $\{v_{i}^{*}\}_{i=0}^{d}$ is diagonal, we have $A^{*}s^{*}=s^{*}A^{*}$.  Recall that the matrix representing $A^{*}$ with respect to the basis $\{v_{i}^{*}\}_{i=0}^{d}$ is irreducible tridiagonal.  From these facts, we have that the matrix representing $A^{*}$ with respect to the basis $\{s^{*}.v_{i}\}_{i=0}^{d}$ is irreducible tridiagonal.  Therefore $\{s^{*}.v_{i}\}_{i=0}^{d}$ is a standard basis for $V$ and $\{\K s^{*}.v_{i}\}_{i=0}^{d}$ is a standard decomposition of $V$.  Recall that $\{\K v_{i}\}_{i=0}^{d}$ and $\{\K v_{d-i}\}_{i=0}^{d}$ are the only standard decompositions of $V$.  Therefore $\{\K s^{*}.v_{i}\}_{i=0}^{d}$ is equal to either $\{\K v_{i}\}_{i=0}^{d}$ or $\{\K v_{d-i}\}_{i=0}^{d}$.  By applying $A$ to bases for $V$ corresponding to each of these decompositions of $V$, we routinely find that the decompositions $\{\K s^{*}.v_{i}\}_{i=0}^{d},\{\K v_{d-i}\}_{i=0}^{d}$ coincide.  It follows that $\theta_{i}=-\theta_{d-i}$ for $0\leq i\leq d$, as desired.

(ii) Similar.
\hfill $\Box$ \\

Recall the Askey-Wilson relations from lines (\ref{E:AW1}), (\ref{E:AW2}).  We now refine these relations in the case where $A,A^{*}$ is totally B/AB.

\begin{theorem}\label{T:BLP}
With reference to Notation \ref{N:notation}, assume $A,A^{*}$ is totally bipartite.  Then $\gamma=\gamma^{*}=\omega=\eta=\eta^{*}$ from Lemma \ref{L:terw-vid} are all zero provided the diameter $d\geq2$.
\end{theorem}

\noindent {\it Proof:} 
Let $\{a_{i}^{*}\}_{i=0}^{d}$ be as in Notation \ref{N:ais}.  By Lemma \ref{L:formulae2} and the fact that $a_{i}^{*}=0$ for $0\leq i\leq d$, the left-hand side of (\ref{E:TV2}) is equal to zero.  Note that the right-hand side of (\ref{E:TV2}) is a quadratic polynomial in $\theta_{i}$.  Then $\theta_{i}$ is a root for $0\leq i\leq d$.  Because $d\geq2$, this polynomial has at least three distinct roots and is therefore zero.  Therefore $\gamma^{*}=\omega=\eta=0$.  By a similar argument using equation (\ref{E:TV3}), we find that $\gamma=\eta^{*}=0$.
\hfill $\Box$ \\

\begin{theorem}\label{T:ABLP}
With reference to Notation \ref{N:notation}, \ref{N:ais}, assume $A,A^{*}$ is totally almost bipartite.  Then $\gamma=\gamma^{*}=\omega=\eta=\eta^{*}$ from Lemma \ref{L:terw-vid} are all zero provided the diameter $d\geq3$.
\end{theorem}

\noindent {\it Proof:} 
Let $\{a_{i}^{*}\}_{i=0}^{d}$ be as in Notation \ref{N:ais}.  Without loss of generality, we assume $a_{d}^{*}\ne0$.  Then $a_{i}^{*}=0$ for $0\leq i\leq d-1$.  By this and Lemma \ref{L:formulae2}, the left-hand side of (\ref{E:TV2}) is equal to zero for $0\leq i\leq d-1$.  Note that the right-hand side of (\ref{E:TV2}) is a quadratic polynomial in $\theta_{i}$.  Then $\theta_{i}$ is a root for $0\leq i\leq d-1$.  Because $d\geq3$, this polynomial has at least three distinct roots and is therefore zero.  Therefore $\gamma^{*}=\omega=\eta=0$.  By a similar argument using equation (\ref{E:TV3}), we find that $\gamma=\eta^{*}=0$.
\hfill $\Box$ \\

We will show that, when $A,A^{*}$ is totally B/AB, the scalars $\varrho,\varrho^{*}$ are nonzero.  In the following Lemma we assume one of $\varrho,\varrho^{*}$ is equal to zero and investigate the consequences.

\begin{lemma}\label{L:rho0}
With reference to Notation \ref{N:notation} and Lemma \ref{L:terw-vid}, the following (i), (ii) hold.
\begin{itemize}
\item[\rm (i)] Suppose the parameters $\gamma,\varrho$ from Lemma \ref{L:terw-vid} are zero.  Then $\{\theta_{i}\}_{i=0}^{d}$ is a geometric progression.  Let $q$ denote the common value of $\theta_{i}/\theta_{i-1}$.  Then $q+q^{-1}=\beta$ from Lemma \ref{L:terw-vid}.

\item[\rm (ii)] Suppose the parameters $\gamma^{*},\varrho^{*}$ from Lemma \ref{L:terw-vid} are zero.  Then $\{\theta^{*}_{i}\}_{i=0}^{d}$ is a geometric progression.  Let $q$ denote the common value of $\theta^{*}_{i}/\theta^{*}_{i-1}$.  Then $q+q^{-1}=\beta$ from Lemma \ref{L:terw-vid}.
\end{itemize}
\end{lemma}

\noindent {\it Proof:} 
(i) Let $r\in\K$ denote a solution to $r+r^{-1}=\beta$.  Substituting $r+r^{-1}$ for $\beta$ in the left-hand equation of (\ref{E:TV1}), and setting $\gamma=\varrho=0$, we find that, for $1\leq i\leq d$,
\begin{align*}
0=&\theta_{i-1}^{2}-(r+r^{-1})\theta_{i-1}\theta_{i}+\theta_{i}^{2}\\
=&(\theta_{i}-r\theta_{i-1})(\theta_{i}-r^{-1}\theta_{i-1}),
\end{align*}
so $\theta_{i}=r\theta_{i-1}$ or $\theta_{i}=r^{-1}\theta_{i-1}$.  Since $\{\theta_{i}\}_{i=0}^{d}$ are mutually distinct, either $\theta_{i}=r\theta_{i-1}$ for all $i$ or $\theta_{i}=r^{-1}\theta_{i-1}$ for $1\leq i\leq d$.  In the former case, set $q=r$ and in the latter case set $q=r^{-1}$.  The result follows.

(ii) Similar.
\hfill $\Box$ \\

\begin{lemma}\label{L:Brho0}
With reference to Notation \ref{N:notation}, assume $A,A^{*}$ is totally bipartite and the diameter $d\geq2$.  Then the scalars $\varrho,\varrho^{*}$ from Lemma \ref{L:terw-vid} are nonzero.
\end{lemma}

\noindent {\it Proof:} 
Assume otherwise.  Without loss of generality, we may assume $\varrho=0$.  By Lemma \ref{L:bipartite}, $\theta_{0}=-\theta_{d}$ and $\theta_{1}=-\theta_{d-1}$.  By Lemma \ref{L:rho0}(i), there exists a nonzero scalar $q$ such that $\theta_{i}=q^{i}\theta_{0}$ for $0\leq i\leq d$ and $q+q^{-1}=\beta$.  Therefore $\theta_{0}=-q^{d}\theta_{0}$ and $q\theta_{0}=-q^{d-1}\theta_{0}$.  From this we obtain $q^{2}\theta_{0}=-q^{d}\theta_{0}=\theta_{0}$.  Therefore $\theta_{0}=\theta_{2}$, a contradiction.  The result follows.
\hfill $\Box$ \\

\begin{lemma}\label{L:special}
With reference to Notation \ref{N:notation}, assume $A,A^{*}$ is totally almost bipartite and $d\geq3$.  Then at least one of  $P(\theta_{0},\theta_{0}),P(\theta_{d},\theta_{d})$ is zero and at least one of $P^{*}(\theta^{*}_{0},\theta^{*}_{0}),P^{*}(\theta^{*}_{d},\theta^{*}_{d})$ is zero.
\end{lemma}

\noindent {\it Proof:} 
By Theorem \ref{T:ABLP}, the right-hand sides of equations  (\ref{E:TV2}), (\ref{E:TV3}) equal zero for $0\leq i\leq d$.  By construction, one of $a_{0},a_{d}$ is nonzero.  If $a_{0}\ne0$ then $P(\theta_{0},\theta_{0})=0$ and if $a_{d}\ne0$ then $P(\theta_{d},\theta_{d})=0$.  Similarly, one of $a^{*}_{0},a^{*}_{d}$ is nonzero.  If $a^{*}_{0}\ne0$ then $P^{*}(\theta^{*}_{0},\theta^{*}_{0})=0$ and if $a^{*}_{d}\ne0$ then $P^{*}(\theta^{*}_{d},\theta^{*}_{d})=0$.
\hfill $\Box$ \\

\begin{lemma}\label{L:ABrho0}
With reference to Notation \ref{N:notation}, assume $A,A^{*}$ is totally almost bipartite and the diameter $d\geq3$.  Then the scalars $\varrho,\varrho^{*}$ from Lemma \ref{L:terw-vid} are nonzero.
\end{lemma}

\noindent {\it Proof:} 
Assume otherwise.  Without loss of generality, we may assume $\varrho=0$.  By Lemma \ref{L:special}, one of $P(\theta_{0},\theta_{0}),P(\theta_{d},\theta_{d})$ is zero.  Reversing the order of the eigenvalues as necessary, we may assume, without loss of generality, that $P(\theta_{d},\theta_{d})=0$.  By the left-hand equation of (\ref{E:TV1}), we have $(2-\beta)\theta_{d}^{2}=0$.  Therefore, either $\beta=2$ or $\theta_{d}=0$.  By Lemma \ref{L:rho0}(i) and the fact that the $\{\theta_{i}\}_{i=0}^{d}$ are distinct, we have $d\leq1$, a contradiction.  The result follows.
\hfill $\Box$ \\

\begin{theorem}\label{T:refinedAW}
Let $A,A^{*}$ denote a totally B/AB Leonard pair.  Then there exists a sequence of scalars $\beta,\varrho,\varrho^{*}$ in $\K$ with $\varrho,\varrho^{*}$ nonzero such that both
\begin{align}
A^{2}A^{*}-\beta AA^{*}A+A^{*}A^{2}&=\varrho A^{*},\label{E:rAW1}\\
A^{*2}A-\beta A^{*}AA^{*}+AA^{*2}&=\varrho^{*}A.\label{E:rAW2}
\end{align}
\end{theorem}

\noindent {\it Proof:} 
Note that equations (\ref{E:rAW1}), (\ref{E:rAW2}) are what we get upon setting $\gamma,\gamma^{*},\omega,\eta,\eta^{*}$ equal to zero in equations (\ref{E:AW1}), (\ref{E:AW2}).  If $d\geq3$, we know that (\ref{E:rAW1}), (\ref{E:rAW2}) hold by Theorems \ref{T:BLP}, \ref{T:ABLP} and Lemmas \ref{L:Brho0}, \ref{L:ABrho0}.  If $d\leq2$ we routinely verify the assertion using Lemmas \ref{L:formulae2} and \ref{L:special}.
\hfill $\Box$ \\

In \cite[Example 5.14]{terwilliger2} a Leonard pair is said to be of \emph{Bannai/Ito type} whenever the common value of (\ref{E:3trec}) is equal to $-1$.  When this occurs, the parameter $\beta$ from Lemma \ref{L:terw-vid} is equal to $-2$ and the relations (\ref{E:rAW1}), (\ref{E:rAW2}) become
\begin{align}
A^{2}A^{*}+2AA^{*}A+A^{*}A^{2}&=\varrho A^{*},\label{E:BIAW1}\\
A^{*2}A+2A^{*}AA^{*}+AA^{*2}&=\varrho^{*}A.\label{E:BIAW2}
\end{align}
When $\varrho,\varrho^{*}$ are equal to $4$, these are equations (\ref{E:xy1}), (\ref{E:xy2}).  Consequently the Leonard pairs associated with the Leonard triple from Theorem \ref{T:Ltrip} are of Bannai/Ito type.

Note that, for $d\geq3$, $\beta$ is uniquely determined in both Lemma \ref{L:terw-vid} and Theorem \ref{T:refinedAW}.  However, when $d\leq 2$, $\beta$ not unique.  As such it is possible for a totally B/AB Leonard pair $A,A^{*}$ with diameter at most $2$ to satisfy equations (\ref{E:AW1}), (\ref{E:AW2}) with $\beta=-2$, but only satisfy equations (\ref{E:rAW1}), (\ref{E:rAW2}) when $\beta\ne-2$.  Because of this, some of the following theorems assume $d\geq3$.

\begin{theorem}\label{T:modLPs}
Let $V$ denote a finite-dimensional irreducible $\mathcal{A}$-module and let $\xi,\xi^{*}$ in $\K$ be nonzero.  Then $\xi x,\xi^{*}y$ act on $V$ as a Leonard pair of Bannai/Ito type.  If $V$ is of type $B$ then the Leonard pair is totally bipartite.  If $V$ is of type $AB$ then the Leonard pair is totally almost bipartite.
\end{theorem}

\noindent {\it Proof:} 
Immediate.
\hfill $\Box$ \\

\begin{theorem}\label{T:Lpairs}
Let $A,A^{*}$ denote a totally B/AB Leonard pair of Bannai/Ito type with diameter $d\geq3$.  Then there exists an irreducible $\mathcal{A}$-module structure on $V$ and nonzero scalars $\xi,\xi^{*}$ such that $A,A^{*}$ act as $\xi x,\xi^{*}y$ respectively.  If $A,A^{*}$ is totally bipartite then $V$ is of type $B$ and if $A,A^{*}$ is totally almost bipartite then $V$ is of type $AB$.  There exist exactly four choices for the scalars $\xi,\xi^{*}$ and the $\mathcal{A}$-module structure.  The scalars $\xi,\xi^{*}$ are each unique up to sign and the $\mathcal{A}$-module structure is uniquely determined by $\xi,\xi^{*}$.
\end{theorem}

\noindent {\it Proof:} 
Since $\K$ is algebraically closed, there exist scalars $\xi,\xi^{*}$ in $\K$ such that $4\xi^{2}=\varrho$, $4\xi^{*2}=\varrho^{*}$.  Because the scalars $\varrho,\varrho^{*}$ are nonzero, the scalars $\xi,\xi^{*}$ are nonzero.  Let $x,y$ act as $A\xi^{-1},A^{*}\xi^{*-1}$ respectively.  By (\ref{E:BIAW1}), (\ref{E:BIAW2}), we find that $x,y$ satisfy (\ref{E:xy1}), (\ref{E:xy2}).

The proof that $V$ is irreducible as an $\mathcal{A}$-module is similar to the proof of Lemma \ref{L:Bd}.  By Theorems \ref{T:trace}, \ref{T:class}, we find that the $\mathcal{A}$-module $V$ is of type $B$ whenever $A,A^{*}$ is totally bipartite and of type $AB$ whenever $A,A^{*}$ is totally almost bipartite.

Given scalars $\xi,\xi^{*}\in\K$, there is at most one $\mathcal{A}$-module structure on $V$ such that $A,A^{*}$ act as $\xi x,\xi^{*}y$ respectively.  Because $\varrho=4\xi^{2}$ and $\varrho^{*}=4\xi^{*2}$ the choices of $\xi,\xi^{*}$ are each unique up to sign.
\hfill $\Box$ \\

Theorem \ref{T:Lpairs} implies the following result of independent interest.

\begin{corollary}\label{C:BABLT}
Let $A,A^{*}$ denote a totally bipartite (resp. totally almost bipartite) Leonard pair of Bannai/Ito type with diameter at least $3$ on the vector space $V$.  Then there exists a linear transformation $A^{\varepsilon}\in\rm{End}$$(V)$ such that $A,A^{*},A^{\varepsilon}$ is a totally bipartite (resp. totally almost bipartite) Leonard triple.
\end{corollary}

\noindent {\it Proof:} 
Let $\xi,\xi^{*}$ be as in Theorem \ref{T:Lpairs} and let $V$ be given the corresponding $\mathcal{A}$-module structure.  Let $\xi^{\varepsilon}\in\K$ be nonzero and let $A^{\varepsilon}$ act as $\xi^{\varepsilon}z$.  By Corollary \ref{C:Ltrip}, $A,A^{*},A^{\varepsilon}$ is a totally bipartite (resp. totally almost bipartite) Leonard triple as desired.
\hfill $\Box$ \\

We now classify the totally B/AB Leonard pairs of Bannai/Ito type with diameter $d\geq3$.  We will be using the notion of isomorphism of Leonard pairs.  For a precise definition, see \cite[Definition 3.4]{terwilligerq}.

\begin{theorem}\label{T:BLpairs}
Let $d$ denote an integer at least $3$ and let $\varrho,\varrho^{*}$ denote scalars in $\K$.  Then the following (i), (ii) are equivalent.
\begin{itemize}
\item[\rm (i)] There exists a totally bipartite Leonard pair $A,A^{*}$ of Bannai/Ito type with diameter $d$ that satisfies equations (\ref{E:BIAW1}), (\ref{E:BIAW2}).

\item[\rm (ii)] The integer $d$ is even and the scalars $\varrho,\varrho^{*}$ are nonzero.
\end{itemize}
Moreover, assume  (i), (ii) hold.  Then the Leonard pair is unique up to isomorphism.
\end{theorem}

\noindent {\it Proof:} 
(ii)$\Rightarrow$(i): Let $V$ denote a finite-dimensional irreducible $\mathcal{A}$-module of type $B(d)$.  Let $\xi,\xi^{*}$ in $\K$ satisfy $4\xi^{2}=\varrho$ and $4\xi^{*2}=\varrho^{*}$.  Let $A,A^{*}$ denote the actions on $V$ of $\xi x,\xi^{*}y$ respectively.  Then, by Theorem \ref{T:modLPs}, $A,A^{*}$ is a totally bipartite Leonard pair of Bannai/Ito type with diameter $d$ that satisfies equations (\ref{E:BIAW1}), (\ref{E:BIAW2}).

(i)$\Rightarrow$(ii): Let $V$ denote the vector space underlying $A,A^{*}$.  By Theorem \ref{T:Lpairs},  there exists an $\mathcal{A}$-module structure on $V$ of type $B$ and nonzero scalars $\xi,\xi^{*}$ such that $A,A^{*}$ act as $\xi x,\xi^{*}y$ respectively.  The dimension of $V$ is $d+1$, so $V$ is of type $B(d)$.  By this and Theorem \ref{T:class}, $d$ is even.  We routinely find that $\varrho=4\xi^{2}$ and $\varrho^{*}=4\xi^{*2}$, so $\varrho,\varrho^{*}$ are nonzero.

Now assume (i), (ii) hold. We show the Leonard pair $A,A^{*}$ is unique up to isomorphism.  Let $B,B^{*}$ denote a totally bipartite Leonard pair of Bannai/Ito type with diameter $d$ that satisfies equations (\ref{E:BIAW1}), (\ref{E:BIAW2}).  We show the Leonard pairs $A,A^{*}$ and $B,B^{*}$ are isomorphic.  Let $V$ denote the vector space underlying $A,A^{*}$ and let $W$ denote the vector space underlying $B,B^{*}$.  By Theorem  \ref{T:Lpairs}, there exist scalars $\xi,\xi^{*}$ in $\K$ and an $\mathcal{A}$-module structure on $V$ such that $A,A^{*}$ act on $V$ as $\xi x,\xi^{*}y$ respectively.  Similarly, there exist scalars $\xi',\xi^{*\prime}$ in $\K$ and an $\mathcal{A}$-module structure on $W$ such that $B,B^{*}$ act on $W$ as $\xi' x,\xi^{*\prime}y$ respectively.  The $\mathcal{A}$-modules $V,W$ are both of type $B(d)$ and hence isomorphic.  By Theorem \ref{T:Lpairs} the scalars $\xi,\xi^{*}$ are unique up to sign, as are the scalars $\xi',\xi^{*\prime}$.  Moreover, both $4\xi^{2},4\xi^{\prime2}$ are equal to $\varrho$ and both $4\xi^{*2},4\xi^{*\prime2}$ are equal to $\varrho^{*}$.  Changing the signs of $\xi,\xi^{*}$ as necessary, we may assume, without loss of generality, that $\xi=\xi'$ and $\xi^{*}=\xi^{*\prime}$.  Let $\phi:V\to W$ denote an isomorphism of $\mathcal{A}$-modules.  Then $\phi\circ A=\xi(\phi\circ x)=\xi(x\circ\phi)=B\circ\phi$ and $\phi\circ A^{*}=\xi^{*}(\phi\circ y)=\xi^{*}(y\circ\phi)=B^{*}\circ\phi$ on $V$.  These equations show the Leonard pairs $A,A^{*}$ and $B,B^{*}$ are isomorphic.
\hfill $\Box$ \\

\begin{theorem}\label{T:ABLpairs}
Let $d$ denote an integer at least $3$ and let $\tau,\tau^{*}$ denote scalars in $\K$.  Then the following (i), (ii) are equivalent.
\begin{itemize}
\item[\rm (i)] There exists a totally almost bipartite Leonard pair $A,A^{*}$ of Bannai/Ito type with diameter $d$, $\mathrm{tr}(A)=\tau$ and $\mathrm{tr}(A^{*})=\tau^{*}$.

\item[\rm (ii)] The scalars $\tau,\tau^{*}$ are nonzero.
\end{itemize}
Moreover, assume  (i), (ii) hold.  Then the Leonard pair is unique up to isomorphism.
\end{theorem}

\noindent {\it Proof:} 
(ii)$\Rightarrow$(i): Let $V$ denote a finite-dimensional irreducible $\mathcal{A}$-module of type $AB(d,0)$.  Let $A,A^{*}$ denote the actions on $V$ of $\tau(-1)^{d}(d+1)^{-1}x,\tau^{*}(-1)^{d}(d+1)^{-1}y$ respectively.  Then, by Theorem \ref{T:modLPs}, $A,A^{*}$ is a totally almost bipartite Leonard pair of Bannai/Ito type with diameter $d$.  By Theorem \ref{T:trace}, $\mathrm{tr}(A)=\tau$ and $\mathrm{tr}(A^{*})=\tau^{*}$.

(i)$\Rightarrow$(ii): Immediate from Definition \ref{D:bipartiteLT}.

Now assume (i), (ii) hold.  We show the Leonard pair $A,A^{*}$ is unique up to isomorphism.  Let $B,B^{*}$ denote a totally almost bipartite Leonard pair of Bannai/Ito type with diameter $d$ such that $\mathrm{tr}(B)=\tau$ and $\mathrm{tr}(B^{*})=\tau^{*}$.  We show the Leonard pairs $A,A^{*}$ and $B,B^{*}$ are isomorphic.  Let $V$ denote the vector space underlying $A,A^{*}$ and let $W$ denote the vector space underlying $B,B^{*}$.  By Theorem  \ref{T:Lpairs}, there exist scalars $\xi,\xi^{*}$ in $\K$ and an $\mathcal{A}$-module structure on $V$ such that $A,A^{*}$ act on $V$ as $\xi x,\xi^{*}y$ respectively.  Similarly, there exist scalars $\xi',\xi^{*\prime}$ in $\K$ and an $\mathcal{A}$-module structure on $W$ such that $B,B^{*}$ act on $W$ as $\xi' x,\xi^{*\prime}y$ respectively.  The $\mathcal{A}$-module $V$ is of type $AB(d,n)$  and the $\mathcal{A}$-module $W$ is of type $AB(d,n')$ for some $n,n'\in\I$.  By Theorem \ref{T:trace} together with $\mathrm{tr}(A)=\mathrm{tr}(B)$ and $\mathrm{tr}(A^{*})=\mathrm{tr}(B^{*})$, we obtain $\xi=\pm\xi'$ and $\xi^{*}=\pm\xi^{*\prime}$, with equality if and only if $n=n'$.  By Theorem \ref{T:Lpairs}, our choice of scalars $\xi,\xi^{*}$ was unique up to sign.  Changing the signs of $\xi,\xi^{*}$ as necessary, we may assume, without loss of generality, that $\xi=\xi'$, $\xi^{*}=\xi^{*\prime}$, and hence $n=n'$.  Then the $\mathcal{A}$-modules $V$ and $W$ are isomorphic.  Let $\phi:V\to W$ denote an isomorphism of $\mathcal{A}$-modules.  Then $\phi\circ A=\xi(\phi\circ x)=\xi(x\circ\phi)=B\circ\phi$ and $\phi\circ A^{*}=\xi^{*}(\phi\circ y)=\xi^{*}(y\circ\phi)=B^{*}\circ\phi$ on $V$.  These equations show the Leonard pairs $A,A^{*}$ and $B,B^{*}$ are isomorphic.
\hfill $\Box$ \\

Note that, given a totally almost bipartite Leonard pair $A,A^{*}$ of Bannai/Ito type with $\tau,\tau^{*}$ from Theorem \ref{T:ABLpairs} and $\varrho,\varrho^{*}$ from equations (\ref{E:BIAW1}), (\ref{E:BIAW2}), we find that
\[
\varrho=\frac{4\tau^{2}}{(d+1)^{2}},\qquad\qquad\varrho^{*}=\frac{4\tau^{*2}}{(d+1)^{2}}.
\]
Given an integer $d$ at least $3$ and nonzero scalars $\varrho,\varrho^{*}$, the scalars $\tau,\tau^{*}$ that satisfy the above equation are each unique up to sign.  Therefore, for each sequence $d,\varrho,\varrho^{*}$ with $d$ an integer at least $3$ and $\varrho,\varrho^{*}$ nonzero, there are exactly $4$ isomorphism classes of totally almost bipartite Leonard pairs of Bannai/Ito type with diameter $d$ that satisfy equations (\ref{E:BIAW1}), (\ref{E:BIAW2}).  Moreover, given a totally bipartite Leonard pair $A,A^{*}$ of Bannai/Ito type, by Definition \ref{D:bipartiteLP}, we find that $\mathrm{tr}(A)=0$ and $\mathrm{tr}(A^{*})=0$.

In Theorem \ref{T:Ltrips} we display a correspondence between totally B/AB Leonard triples of Bannai/Ito type and $\mathcal{A}$-modules.  To do this, we present further results about Leonard pairs.  Let $A,A^{*}$ denote a Leonard pair.  With reference to Notation \ref{N:notation}, let $E_{i},E_{i}^{*}$ denote the primitive idempotents corresponding to $\theta_{i},\theta_{i}^{*}$ respectively for $0\leq i\leq d$.

\begin{lemma}\label{L:lin-ind}
With reference to Notation \ref{N:notation}, let $A,A^{*}$ be totally B/AB and of Bannai/Ito type with diameter $d\geq3$.  Then the elements
\begin{equation}\label{E:brown}
A^{2}A^{*},\qquad\qquad AA^{*}A,\qquad\qquad A^{*}A^{2},
\end{equation}
are linearly independent.
\end{lemma}

\noindent {\it Proof:} 
Let $s,t,u$  be scalars in $\K$ satisfying $sA^{2}A^{*}+tAA^{*}A+uA^{*}A^{2}=0$.  We show that each of $s,t,u$ is zero.  The following hold:
\begin{align}
sE_{0}^{*}A^{2}A^{*}E_{0}^{*}+tE_{0}^{*}AA^{*}AE_{0}^{*}+uE_{0}^{*}A^{*}A^{2}E_{0}^{*}=&0,\label{E:stuff1}\\
sE_{0}^{*}A^{2}A^{*}E_{2}^{*}+tE_{0}^{*}AA^{*}AE_{2}^{*}+uE_{0}^{*}A^{*}A^{2}E_{2}^{*}=&0,\label{E:stuff2}\\
sE_{2}^{*}A^{2}A^{*}E_{0}^{*}+tE_{2}^{*}AA^{*}AE_{0}^{*}+uE_{2}^{*}A^{*}A^{2}E_{0}^{*}=&0.\label{E:stuff3}
\end{align}

With respect to the basis $\{v_{i}^{*}\}_{i=0}^{d}$ from Notation \ref{N:notation}, the matrix representing $A^{*}$ is diagonal with $(i,i)$-entry $\theta_{i}^{*}$ for $0\leq i\leq d$.  For $0\leq i\leq d$, the matrix representing $E_{i}^{*}$ has $(i,i)$-entry $1$ and all other entries zero.  With respect to Notation \ref{N:ais}, the matrix representing $A$ is irreducible tridiagonal with $(i,i)$ entry $a_{i}$ for $0\leq i\leq d$.  The matrix representing $A$ is either bipartite or almost bipartite.  Therefore at most one of $a_{0},a_{d}$ is nonzero and $a_{i}=0$ for $1\leq i\leq d-1$.  Reversing the order of the eigenvalues as necessary we may assume, without loss of generality, that $a_{0}=0$.  Based on this information, we routinely find that equations (\ref{E:stuff1})--(\ref{E:stuff3}) reduce to

\begin{align}
(s\theta_{0}^{*}+t\theta_{1}^{*}+u\theta_{0}^{*})E_{0}^{*}A^{2}E_{0}^{*}=&0,\\
(s\theta_{2}^{*}+t\theta_{1}^{*}+u\theta_{0}^{*})E_{0}^{*}A^{2}E_{2}^{*}=&0,\\
(s\theta_{0}^{*}+t\theta_{1}^{*}+u\theta_{2}^{*})E_{2}^{*}A^{2}E_{0}^{*}=&0.
\end{align}

Moreover, $E_{0}^{*}A^{2}E_{0}^{*},E_{0}^{*}A^{2}E_{2}^{*},E_{2}^{*}A^{2}E_{0}^{*}$ are all nonzero, resulting in the following equations:

\begin{align}
s\theta_{0}^{*}+t\theta_{1}^{*}+u\theta_{0}^{*}=&0,\label{E:system1}\\
s\theta_{2}^{*}+t\theta_{1}^{*}+u\theta_{0}^{*}=&0,\label{E:system2}\\
s\theta_{0}^{*}+t\theta_{1}^{*}+u\theta_{2}^{*}=&0.\label{E:system3}
\end{align}

We view (\ref{E:system1})--(\ref{E:system3}) as a system of linear equations in the indeterminates $s,t,u$.  The determinant of the coefficient matrix is $-\theta_{1}^{*}(\theta_{0}^{*}-\theta_{2}^{*})^{2}$.  Because $\{\theta_{i}^{*}\}_{i=0}^{d}$ are distinct, we have $\theta_{0}^{*}-\theta_{2}^{*}\ne0$.  Combining Theorem \ref{T:Lpairs}, the eigenvalue data for $y$ from Propositions \ref{P:idemtable1}, \ref{P:idemtable2} and the fact that $d\geq3$, we find that $\theta_{1}^{*}\ne0$.  From this, we routinely find that $s=0,t=0,u=0$ is the only solution to the system (\ref{E:system1})--(\ref{E:system3}).  Therefore, (\ref{E:brown}) are linearly independent as desired.
\hfill $\Box$ \\

Now, let $\mathcal{X}$ denote the $\K$-subspace of $V$ consisting of the $X\in\mathrm{End}(V)$ such that both
\begin{align}
E_{i}XE_{j}=&0\qquad\qquad\mbox{if $|i-j|>1,$}\\
E_{i}^{*}XE_{j}^{*}=&0\qquad\qquad\mbox{if $|i-j|>1,$}
\end{align}
for $0\leq i,j\leq d$.  Observe that, if $A,A^{*},A^{\varepsilon}$ is a Leonard triple, then $A^{\varepsilon}\in\mathcal{X}$.

\begin{lemma}{\normalfont\cite[Theorem 1.5]{nom-terw}}\label{L:nom-terw}
The space $\mathcal{X}$ is spanned by
\begin{equation}\label{E:nom-terw}
I, A, A^{*}, AA^{*},A^{*}A.
\end{equation}
Moreover, (\ref{E:nom-terw}) is a basis for $\mathcal{X}$ provided $d\geq2$.
\end{lemma}

\section{Totally B/AB Leonard triples of Bannai/Ito type}\label{S:LTclassification}

In Section \ref{S:LPclassification}, we classified the Leonard pairs arising from finite-dimensional irreducible $\mathcal{A}$-modules.  In this section, we classify the Leonard triples arising from finite-dimensional irreducible $\mathcal{A}$-modules.  We show that they correspond to a family of totally B/AB Leonard triples said to have Bannai/Ito type.  From this correspondence we classify the totally B/AB Leonard triples of Bannai/Ito type with diameter at least $3$.

\begin{notation}\label{N:LTnotation}
Let $V$ denote a vector space over $\K$ with finite positive dimension.  Let $A,A^{*},A^{\varepsilon}$ denote a Leonard triple on $V$.  Let $\{v_{i}\}_{i=0}^{d}$ denote a basis for $V$ under which $A$ is diagonal and $A^{*},A^{\varepsilon}$ are irreducible tridiagonal. Let $\{v^{*}_{i}\}_{i=0}^{d}$ denote a basis for $V$ under which $A^{*}$ is diagonal and $A^{\varepsilon},A$ are irreducible tridiagonal.  Let $\{v^{\varepsilon}_{i}\}_{i=0}^{d}$ denote a basis for $V$ under which $A^{\varepsilon}$ is diagonal and $A,A^{*}$ are irreducible tridiagonal.  For $0\leq i\leq d$, let $\theta_{i}$ denote the eigenvalue for $A$ associated with $v_{i}$, let $\theta_{i}^{*}$ denote the eigenvalue for $A^{*}$ associated with $v_{i}^{*}$ and let $\theta_{i}^{\varepsilon}$ denote the eigenvalue for $A^{\varepsilon}$ associated with $v_{i}^{\varepsilon}$.
\end{notation}

\begin{definition}\label{D:BILT}
{\normalfont We say that a Leonard triple $A,A^{*},A^{\varepsilon}$ is of \emph{Bannai/Ito type} whenever all of the associated Leonard pairs are of Bannai/Ito type.}
\end{definition}

\begin{lemma}\label{L:BILT}
Let $A,A^{*},A^{\varepsilon}$ denote a Leonard triple.  If any of the six Leonard pairs associated with $A,A^{*},A^{\varepsilon}$ is of Bannai/Ito type, then the Leonard triple $A,A^{*},A^{\varepsilon}$ is of Bannai/Ito type.
\end{lemma}

\noindent {\it Proof:} 
Assume otherwise.  If the Leonard pair $A,A^{*}$ is of Bannai/Ito type then so is $A^{*},A$.  Therefore we may assume, without loss of generality, that the Leonard pair $A,A^{*}$ is of Bannai/Ito type and the Leonard pair $A,A^{\varepsilon}$ is not of Bannai/Ito type.  With reference to Notation \ref{N:LTnotation}, consider the common value of $(\theta_{i-2}-\theta_{i+1})/(\theta_{i-1}-\theta_{i})$ for $2\leq i\leq d-1$.  Because $A,A^{*}$ is of Bannai/Ito type, that common value is equal to $-1$.  Because $A,A^{\varepsilon}$ is not of Bannai/Ito type, that same common value is not equal to $-1$.  This is a contradiction, and the result follows.
\hfill $\Box$ \\

\begin{theorem}\label{T:modLTs}
Let $V$ denote a finite-dimensional irreducible $\mathcal{A}$-module and let $A,A^{*},A^{\varepsilon}$ denote the Leonard triple from Corollary \ref{C:Ltrip}.  Then $A,A^{*},A^{\varepsilon}$ is of Bannai/Ito type.  If $V$ is of type $B$ then $A,A^{*},A^{\varepsilon}$ is totally bipartite.  If $V$ is of type $AB$ then $A,A^{*},A^{\varepsilon}$ is totally almost bipartite.
\end{theorem}

\noindent {\it Proof:} 
Immediate.
\hfill $\Box$ \\

\begin{theorem}\label{T:Ltrips}
Let $A,A^{*},A^{\varepsilon}$ denote a totally B/AB Leonard triple of Bannai/Ito type with diameter $d\geq3$.  Then there exists an irreducible $\mathcal{A}-module$ structure on $V$ and nonzero scalars $\xi,\xi^{*},\xi^{\varepsilon}$ in $\K$ such that $A,A^{*}A^{\varepsilon}$ act as $\xi x,\xi^{*}y,\xi^{\varepsilon}z$ respectively.  If $A,A^{*},A^{\varepsilon}$ is totally bipartite then $V$ is of type $B$ and if $A,A^{*},A^{\varepsilon}$ is totally almost bipartite then $V$ is of type $AB$.  There exist exactly four choices for the scalars $\xi,\xi^{*},\xi^{\varepsilon}$ and the $\mathcal{A}$-module structure.  The scalars $\xi,\xi^{*}$ are each unique up to sign and the scalar $\xi^{\varepsilon}$ and the $\mathcal{A}$-module structure are uniquely determined by $\xi,\xi^{*}$.
\end{theorem}

\noindent {\it Proof:} 
We first claim there exist scalars $\zeta_{1},\zeta_{2},\zeta^{*}_{1},\zeta^{*}_{2},\zeta^{\varepsilon}_{1},\zeta^{\varepsilon}_{2}\in\K$ such that
\begin{align}
\zeta_{1}^{\varepsilon}AA^{*}+\zeta_{2}^{\varepsilon}A^{*}A&=A^{\varepsilon}\label{E:zeta1},\\
\zeta_{1}A^{*}A^{\varepsilon}+\zeta_{2}A^{\varepsilon}A^{*}&=A\label{E:zeta2},\\
\zeta_{1}^{*}A^{\varepsilon}A+\zeta_{2}^{*}AA^{\varepsilon}&=A^{*}\label{E:zeta3}.
\end{align}

To prove the claim, we first show that line (\ref{E:zeta1}) holds.  By Lemma \ref{L:nom-terw} and the fact that $d\geq3$, there exist unique scalars $\alpha_{1},\alpha_{2},\alpha_{3},\alpha_{4},\alpha_{5}\in\K$ such that
\begin{equation}\label{E:zetas}
A^{\varepsilon}=\alpha_{1}I+\alpha_{2}A+\alpha_{3}A^{*}+\alpha_{4}AA^{*}+\alpha_{5}A^{*}A.
\end{equation}

In equation (\ref{E:zetas}) set $\zeta_{1}^{\varepsilon}=\alpha_{4}$ and $\zeta_{2}^{\varepsilon}=\alpha_{5}$.  We show that $\alpha_{1},\alpha_{2},\alpha_{3}$ are equal to zero.  We first show that $\alpha_{1},\alpha_{2}$ are equal to zero.  Consider the matrix $B^{\varepsilon}$ representing $A^{\varepsilon}$ with respect to the basis $\{v_{i}\}_{i=0}^{d}$ from Notation \ref{N:LTnotation}.  By construction, $B^{\varepsilon}$ is irreducible tridiagonal and either bipartite or almost bipartite.  The matrices representing $A^{*},AA^{*},A^{*}A$ are also tridiagonal and either bipartite or almost bipartite.  Therefore, $B^{\varepsilon}_{i,i}=0$ for $1\leq i\leq d-1$.  By (\ref{E:zetas}), we have that, for $1\leq i\leq d-1$,
\[
B^{\varepsilon}_{i,i}=\alpha_{1}+\alpha_{2}\theta_{i}.
\]
Because $d\geq3$ and $\{\theta_{i}\}_{i=1}^{d-1}$ are distinct, we have that $\alpha_{1},\alpha_{2}$ are both equal to zero.  The proof that $\alpha_{3}=0$ is similar, using the matrix representing $A^{\varepsilon}$ with respect to the basis $\{v_{i}^{*}\}_{i=0}^{d}$ from Notation \ref{N:LTnotation}.  Therefore equation (\ref{E:zeta1}) holds.  Equations (\ref{E:zeta2}), (\ref{E:zeta3}) are similar and the claim follows.

We now refine the relations (\ref{E:zeta1})--(\ref{E:zeta3}).  We claim that there exist nonzero scalars $\zeta,\zeta^{*},\zeta^{\varepsilon}\in\K$ such that
\begin{align}
\zeta^{\varepsilon}(AA^{*}+A^{*}A)&=A^{\varepsilon}\label{E:rzeta1},\\
\zeta(A^{*}A^{\varepsilon}+A^{\varepsilon}A^{*})&=A\label{E:rzeta2},\\
\zeta^{*}(A^{\varepsilon}A+AA^{\varepsilon})&=A^{*}\label{E:rzeta3}.
\end{align}

Substituting the left-hand side of equation (\ref{E:zeta1}) for $A^{\varepsilon}$ in equation (\ref{E:zeta3}), we find that
\begin{equation}\label{E:something}
\zeta_{1}^{\varepsilon}\zeta_{2}^{*}A^{2}A^{*}+(\zeta_{1}^{\varepsilon}\zeta_{1}^{*}+\zeta_{2}^{\varepsilon}\zeta_{2}^{*})AA^{*}A+\zeta_{2}^{\varepsilon}\zeta_{1}^{*}A^{*}A^{2}=A^{*}.
\end{equation}
Equations (\ref{E:BIAW1}) and (\ref{E:something}) both express $A^{*}$ as a linear combination of (\ref{E:brown}).  By Lemma \ref{L:lin-ind}, we have
\begin{align}
\varrho\zeta_{1}^{\varepsilon}\zeta_{2}^{*}&=1\label{E:thing1},\\
\varrho(\zeta_{1}^{\varepsilon}\zeta_{1}^{*}+\zeta_{2}^{\varepsilon}\zeta_{2}^{*})&=2\label{E:thing2},\\
\varrho\zeta_{2}^{\varepsilon}\zeta_{1}^{*}&=1\label{E:thing3}.
\end{align}
By equation (\ref{E:thing1}), we have $\zeta_{1}^{\varepsilon}\ne0$ and by equation (\ref{E:thing3}), we have $\zeta_{2}^{\varepsilon}\ne0$.  Solving equations (\ref{E:thing1}) and (\ref{E:thing3}) for $\zeta_{2}^{*},\zeta_{1}^{*}$ respectively and substituting into equation (\ref{E:thing2}), we get $\zeta_{1}^{\varepsilon}(\zeta_{2}^{\varepsilon})^{-1}+\zeta_{2}^{\varepsilon}(\zeta_{1}^{\varepsilon})^{-1}=2$.  Therefore $\zeta_{1}^{\varepsilon}=\zeta_{2}^{\varepsilon}$ and both are nonzero.  Let $\zeta^{\varepsilon}$ denote the common value of $\zeta_{1}^{\varepsilon},\zeta_{2}^{\varepsilon}$.  Then equation (\ref{E:rzeta1}) holds.  Equations (\ref{E:rzeta2}), (\ref{E:rzeta3}) are similar and the second claim follows.

Since $\K$ is algebraically closed and $\zeta,\zeta^{*},\zeta^{\varepsilon}$ are nonzero, there exist $\xi,\xi^{*},\xi^{\varepsilon}$ such that $\xi^{2}=(4\zeta^{*}\zeta^{\varepsilon})^{-1}$, $\xi^{*2}=(4\zeta^{\varepsilon}\zeta)^{-1}$ and $\xi^{\varepsilon2}=(4\zeta\zeta^{*})^{-1}$.  The choices for $\xi,\xi^{*},\xi^{\varepsilon}$ are unique up to sign and $\xi\xi^{*}\xi^{\varepsilon}=\pm(8\zeta\zeta^{*}\zeta^{\varepsilon})^{-1}$.  Choose $\xi,\xi^{*}\xi^{\varepsilon}$ such that $\xi\xi^{*}\xi^{\varepsilon}=(8\zeta\zeta^{*}\zeta^{\varepsilon})^{-1}$.  We have $\xi,\xi^{*},\xi^{\varepsilon}\ne0$.  Let $x,y,z$ act as $A\xi^{-1},A^{*}\xi^{*-1},A^{\varepsilon}\xi^{\varepsilon-1}$ respectively.  By (\ref{E:rzeta1})--(\ref{E:rzeta3}) we have that $x,y,z$ satisfy (\ref{E:rel1})--(\ref{E:rel3}).

The proof that $V$ is irreducible as an $\mathcal{A}$-module is similar to the proof of Lemma \ref{L:Bd}.  By Theorems \ref{T:trace}, \ref{T:class}, we find that, if $A,A^{*},A^{\varepsilon}$ is totally bipartite then $V$ is of type $B$ and if $A,A^{*},A^{\varepsilon}$ is totally almost bipartite then $V$ is of type $AB$.

Given scalars $\xi,\xi^{*},\xi^{\varepsilon}\in\K$, there is at most one $\mathcal{A}$-module structure on $V$ such that $A,A^{*},A^{\varepsilon}$ act as $\xi x,\xi^{*}y,\xi^{\varepsilon}z$ respectively.  Because $\xi^{2}=(4\zeta^{*}\zeta^{\varepsilon})^{-1}$, $\xi^{*2}=(4\zeta^{\varepsilon}\zeta)^{-1}$, $\xi^{\varepsilon2}=(4\zeta\zeta^{*})^{-1}$ and $\xi\xi^{*}\xi^{\varepsilon}=(8\zeta\zeta^{*}\zeta^{\varepsilon})^{-1}$, the choices of $\xi,\xi^{*}$ are unique up to sign change and $\xi^{\varepsilon}$ is uniquely determined by $\xi,\xi^{*}$.
\hfill $\Box$ \\

In Theorem \ref{T:Ltrips} we assume that $d\geq3$.  To see that this assumption is necessary, we show that, for $d=2$, the Theorem is false.  By \cite[Theorems 10.1(i), 10.2(ii), 10.4(iii)]{miklavic} with $d=2$,

\[
A=
\begin{pmatrix}
2 & 0 & 0\\
0 & 0 & 0\\
0 & 0 & -2
\end{pmatrix},\qquad
A^{*}=
\begin{pmatrix}
0 & 2 & 0\\
1 & 0 & 1\\
0 & 2 & 0
\end{pmatrix},\qquad
A^{\varepsilon}=
\begin{pmatrix}
0 & -2\mathbf{i} & 0\\
\mathbf{i} & 0 & -\mathbf{i}\\
0 & 2\mathbf{i} & 0
\end{pmatrix}
\]
is a Leonard triple with diameter $2$.  Observing \cite[Theorems 10.1(ii),(iii), 10.2(i),(iii), 10.4(i), (ii)]{miklavic} we find that the Leonard triple is totally bipartite, and we routinely find that each Leonard pair obtained from this Leonard triple satisfies equations (\ref{E:BIAW1}), (\ref{E:BIAW2}) with $\varrho=4$ and $\varrho^{*}=4$, and is hence of Bannai/Ito type.  However, there are no scalars $\zeta,\zeta^{*},\zeta^{\varepsilon}$ that satisfy equation (\ref{E:rzeta1}).  Therefore, there is no $\mathcal{A}$-module structure as described in Theorem \ref{T:Ltrips}.

We now classify the totally B/AB Leonard triples of Bannai/Ito type with diameter $d\geq3$.  We will be using the notion of isomorphism of Leonard triples.  For a precise definition, see \cite[Definition 8.2]{curtin}.

\begin{theorem}\label{T:BLtriples}
Let $d$ denote an integer at least $3$ and let $\zeta,\zeta^{*},\zeta^{\varepsilon}$ denote scalars in $\K$.  Then the following (i), (ii) are equivalent.
\begin{itemize}
\item[\rm (i)] There exists a totally bipartite Leonard triple $A,A^{*},A^{\varepsilon}$ of Bannai/Ito type with diameter $d$ that satisfies equations (\ref{E:rzeta1})--(\ref{E:rzeta3}).

\item[\rm (ii)] The integer $d$ is even and the scalars $\zeta,\zeta^{*},\zeta^{\varepsilon}$ are nonzero.
\end{itemize}
Moreover, assume  (i), (ii) hold.  Then the Leonard triple is unique up to isomorphism.
\end{theorem}

\noindent {\it Proof:} 
(ii)$\Rightarrow$(i): Let $V$ denote a finite-dimensional irreducible $\mathcal{A}$-module of type $B(d)$.  Let $\xi,\xi^{*},\xi^{\varepsilon}$ in $\K$ satisfy $\xi^{2}=(4\zeta^{*}\zeta^{\varepsilon})^{-1}$, $\xi^{*2}=(4\zeta^{\varepsilon}\zeta)^{-1}$, $\xi^{\varepsilon2}=(4\zeta\zeta^{*})^{-1}$ and $\xi\xi^{*}\xi^{\varepsilon}=(8\zeta\zeta^{*}\zeta^{\varepsilon})^{-1}$.  Let $A,A^{*},A^{\varepsilon}$  denote the actions on $V$ of $\xi x,\xi^{*}y,\xi^{\varepsilon}z$ respectively.  Then, by Theorem \ref{T:modLTs}, $A,A^{*},A^{\varepsilon}$ is a totally bipartite Leonard triple of Bannai/Ito type with diameter $d$ that satisfies equations (\ref{E:rzeta1})--(\ref{E:rzeta3}).

(i)$\Rightarrow$(ii): By Theorem \ref{T:Ltrips},  there exists an $\mathcal{A}$-module structure on $V$ of type $B$ and nonzero scalars $\xi,\xi^{*},\xi^{\varepsilon}$ such that $A,A^{*},A^{\varepsilon}$ act as $\xi x,\xi^{*}y,\xi^{\varepsilon}z$ respectively.  The dimension of $V$ is $d+1$, so $V$ is of type $B(d)$.  By this and Theorem \ref{T:class}, $d$ is even.  We routinely find that $\zeta=\xi(2\xi^{*}\xi^{\varepsilon})^{-1}$, $\zeta^{*}=\xi^{*}(2\xi^{\varepsilon}\xi)^{-1}$ and $\zeta^{\varepsilon}=\xi^{\varepsilon}(2\xi\xi^{*})^{-1}$ so $\zeta,\zeta^{*},\zeta^{\varepsilon}$ are nonzero.

Now assume (i), (ii) hold.  We show the Leonard triple $A,A^{*},A^{\varepsilon}$ is unique up to isomorphism.  Let $B,B^{*},B^{\varepsilon}$ denote a totally bipartite Leonard triple of Bannai/Ito type with diameter $d$ that satisfies equations (\ref{E:rzeta1})--(\ref{E:rzeta3}).  We show the Leonard triples $A,A^{*},A^{\varepsilon}$ and $B,B^{*},B^{\varepsilon}$ are isomorphic.  Let $V$ denote the vector space underlying $A,A^{*},A^{\varepsilon}$ and let $W$ denote the vector space underlying $B,B^{*},B^{\varepsilon}$.  By Theorem  \ref{T:Ltrips}, there exist scalars $\xi,\xi^{*},\xi^{\varepsilon}$ in $\K$ and an $\mathcal{A}$-module structure on $V$ such that $A,A^{*},A^{\varepsilon}$ act on $V$ as $\xi x,\xi^{*}y,\xi^{\varepsilon}z$ respectively.  Similarly, there exist scalars $\xi',\xi^{*\prime},\xi^{\varepsilon\prime}$ in $\K$ and an $\mathcal{A}$-module structure on $W$ such that $B,B^{*},B^{\varepsilon}$ act on $W$ as $\xi'x,\xi^{*\prime}y,\xi^{\varepsilon\prime}z$ respectively.  The $\mathcal{A}$-modules $V,W$ are both of type $B(d)$ and hence isomorphic.  By Theorem 8.16 the scalars $\xi,\xi^{*}$ are unique up to sign as are the scalars $\xi',\xi^{*\prime}$.  Moreover, the scalar $\xi^{\varepsilon}$ is uniquely determined by $\xi,\xi^{*}$ and the scalar $\xi^{\varepsilon\prime}$ is uniquely determined by $\xi',\xi^{*\prime}$.  Moreover, both $\xi^{2},\xi^{\prime2}$ are equal to $(4\zeta^{*}\zeta^{\varepsilon})^{-1}$, both $\xi^{*2},\xi^{*\prime2}$ are equal to $(4\zeta^{\varepsilon}\zeta)^{-1}$,  both $\xi^{\varepsilon2},\xi^{\varepsilon\prime2}$ are equal to $(4\zeta\zeta^{*})^{-1}$ and both $\xi\xi^{*}\xi^{\varepsilon},\xi'\xi^{*\prime}\xi^{\varepsilon\prime}$ are equal to $(8\zeta\zeta^{*}\zeta^{\varepsilon})^{-1}$.  Changing the signs of $\xi,\xi^{*},\xi^{\varepsilon}$ as necessary, we may assume, without loss of generality, that $\xi=\xi'$, $\xi^{*}=\xi^{*\prime}$ and $\xi^{\varepsilon}=\xi^{\varepsilon\prime}$.  Let $\phi:V\to W$ denote an isomorphism of $\mathcal{A}$-modules.  Then $\phi\circ A=\xi(\phi\circ x)=\xi(x\circ\phi)=B\circ\phi$, $\phi\circ A^{*}=\xi^{*}(\phi\circ y)=\xi^{*}(y\circ\phi)=B^{*}\circ\phi$ and $\phi\circ A^{\varepsilon}=\xi^{\varepsilon}(\phi\circ z)=\xi^{\varepsilon}(z\circ\phi)=B^{\varepsilon}\circ\phi$ on $V$.  These equations show the Leonard triples $A,A^{*},A^{\varepsilon}$ and $B,B^{*},B^{\varepsilon}$ are isomorphic.
\hfill $\Box$ \\

\begin{theorem}\label{T:ABLtriples}
Let $d$ denote an integer at least $3$ and let $\tau,\tau^{*},\tau^{\varepsilon}$ denote scalars in $\K$.  Then the following (i), (ii) are equivalent.
\begin{itemize}
\item[\rm (i)] There exists a totally almost bipartite Leonard triple $A,A^{*},A^{\varepsilon}$ of Bannai/Ito type with diameter $d$, $\mathrm{tr}(A)=\tau$, $\mathrm{tr}(A^{*})=\tau^{*}$ and $\mathrm{tr}(A^{\varepsilon})=\tau^{\varepsilon}$.

\item[\rm (ii)] The scalars $\tau,\tau^{*},\tau^{\varepsilon}$ are nonzero.
\end{itemize}
Moreover, assume  (i), (ii) hold.  Then the Leonard triple is unique up to isomorphism.
\end{theorem}

\noindent {\it Proof:} 
(ii)$\Rightarrow$(i): Let $V$ denote a finite-dimensional irreducible $\mathcal{A}$-module of type $AB(d,0)$.  Let $A,A^{*},A^{\varepsilon}$ denote the actions of $\tau(-1)^{d}(d+1)^{-1}x,\tau^{*}(-1)^{d}(d+1)^{-1}y,A^{\varepsilon}=\tau^{\varepsilon}(-1)^{d}(d+1)^{-1}z$ respectively.  Then, by Theorem \ref{T:modLPs}, $A,A^{*},A^{\varepsilon}$ is a totally almost bipartite Leonard triple of Bannai/Ito type with diameter $d$.  By Theorem \ref{T:trace}, $\mathrm{tr}(A)=\tau$, $\mathrm{tr}(A^{*})=\tau^{*}$ and $\mathrm{tr}(A^{\varepsilon})=\tau^{\varepsilon}$.

(i)$\Rightarrow$(ii): Immediate from Definition \ref{D:bipartiteLT}.

Now assume (i), (ii) hold.  We show the Leonard triple $A,A^{*},A^{\varepsilon}$ is unique up to isomorphism.  Let $B,B^{*},B^{\varepsilon}$ denote a totally almost bipartite Leonard triple of Bannai/Ito type with diameter $d$ such that $\mathrm{tr}(B)=\tau$, $\mathrm{tr}(B^{*})=\tau^{*}$ and $\mathrm{tr}(B^{\varepsilon})=\tau^{\varepsilon}$.  We show the Leonard triples $A,A^{*},A^{\varepsilon}$ and $B,B^{*},B^{\varepsilon}$ are isomorphic.  Let $V$ denote the vector space underlying $A,A^{*},A^{\varepsilon}$ and let $W$ denote the vector space underlying $B,B^{*},B^{\varepsilon}$.  By Theorem  \ref{T:Ltrips}, there exist scalars $\xi,\xi^{*},\xi^{\varepsilon}$ in $\K$ and an $\mathcal{A}$-module structure on $V$ such that $A,A^{*},A^{\varepsilon}$ act on $V$ as $\xi x,\xi^{*}y,\xi^{\varepsilon}z$ respectively.  Similarly, there exist scalars $\xi',\xi^{*\prime},\xi^{\varepsilon\prime}$ in $\K$ and an $\mathcal{A}$-module structure on $W$ such that $B,B^{*},B^{\varepsilon}$ act on $W$ as $\xi' x,\xi^{*\prime}y,\xi^{\varepsilon\prime}z$ respectively.  The $\mathcal{A}$-module $V$ is of type $AB(d,n)$  and the $\mathcal{A}$-module $W$ is of type $AB(d,n')$ for some $n,n'\in\I$.  By Theorem \ref{T:trace} together with $\mathrm{tr}(A)=\mathrm{tr}(B)$, $\mathrm{tr}(A^{*})=\mathrm{tr}(B^{*})$ and $\mathrm{tr}(A^{\varepsilon})=\mathrm{tr}(B^{\varepsilon})$, we obtain $\xi=\pm\xi'$, $\xi^{*}=\pm\xi^{*\prime}$ and $\xi^{\varepsilon}=\pm\xi^{\varepsilon\prime}$, with equality if and only if $n=n'$.  By Theorem \ref{T:Ltrips}, our choice of scalars $\xi,\xi^{*}$ was unique up to sign and our choice of $\xi^{\varepsilon}$ was determined by $\xi,\xi^{*}$.  Changing the signs of $\xi,\xi^{*},\xi^{\varepsilon}$ as necessary, we may assume, without loss of generality, that $\xi=\xi'$, $\xi^{*}=\xi^{*\prime}$, $\xi^{\varepsilon}=\xi^{\varepsilon\prime}$ and hence $n=n'$.  Then the $\mathcal{A}$-modules $V$ and $W$ are isomorphic.  Let $\phi:V\to W$ denote an isomorphism of $\mathcal{A}$-modules.  Then $\phi\circ A=\xi(\phi\circ x)=\xi(x\circ\phi)=B\circ\phi$, $\phi\circ A^{*}=\xi^{*}(\phi\circ y)=\xi^{*}(y\circ\phi)=B^{*}\circ\phi$ and $\phi\circ A^{\varepsilon}=\xi^{\varepsilon}(\phi\circ z)=\xi^{\varepsilon}(z\circ\phi)=B^{\varepsilon}\circ\phi$ on $V$.  These equations show the Leonard triples $A,A^{*},A^{\varepsilon}$ and $B,B^{*},B^{\varepsilon}$ are isomorphic.
\hfill $\Box$ \\

Note that, given a totally almost bipartite Leonard triple $A,A^{*},A^{\varepsilon}$ of Bannai/Ito type with $\tau,\tau^{*},\tau^{\varepsilon}$ from Theorem \ref{T:ABLtriples} and $\zeta,\zeta^{*},\zeta^{\varepsilon}$ from equations (\ref{E:rzeta1})--(\ref{E:rzeta3}), we find that
\[
\zeta=\frac{(-1)^{d}(d+1)\tau}{2\tau^{*}\tau^{\varepsilon}},\qquad\zeta^{*}=\frac{(-1)^{d}(d+1)\tau^{*}}{2\tau^{\varepsilon}\tau},\qquad\zeta^{\varepsilon}=\frac{(-1)^{d}(d+1)\tau^{\varepsilon}}{2\tau\tau^{*}}.
\]
Given an integer $d$ at least three and nonzero scalars $\zeta,\zeta^{*},\zeta^{\varepsilon}$, the scalars $\tau,\tau^{*},\tau^{\varepsilon}$ that satisfy the above equation are unique up to changing the sign of an even number of them.  Therefore, for each sequence $d,\zeta,\zeta^{*},\zeta^{\varepsilon}$ with $d$ an integer at least $3$ and $\zeta,\zeta^{*},\zeta^{\varepsilon}$ nonzero, there are exactly $4$ isomorphism classes of totally almost bipartite Leonard triples of Bannai/Ito type with diameter $d$ that satisfy equations (\ref{E:rzeta1})--(\ref{E:rzeta3}).  Moreover, given a totally bipartite Leonard triple $A,A^{*},A^{\varepsilon}$ of Bannai/Ito type, by Definition \ref{D:bipartiteLT}, $\mathrm{tr}(A)=0$, $\mathrm{tr}(A^{*})=0$ and $\mathrm{tr}(A^{\varepsilon})=0$.

\section{Acknowledgment}\label{S:acknowledgment}

This paper was written while the author was a graduate student at the University of Wisconsin-Madison.  The author would like to thank his advisor, Paul Terwilliger, for offering many valuable ideas and suggestions.

\end{document}